\documentclass{article}

\usepackage{amsfonts}
\usepackage{amssymb}

\usepackage{amsmath}

\begin{document}

\title{\bf Twists in $U\left(\mathfrak{sl}_3\right)$ and their quantizations}

\author{\bf P P Kulish$^1$,\\ \bf V D Lyakhovsky$^2$,\\ \bf M E Samsonov$^3$\\[3mm]
$^1$ St Petersburg Dept. of Steklov Mathematical Institute,\\[2mm]
St Petersburg, 191011, Russia\\[2mm]
$^2$ $^3$ Theoretical Department, St Petersburg State University,\\
St Petersburg, 198904, Russia\\[2mm]
e-mail addresses:\\ kulish@pdmi.ras.ru,\\ lyakhovs@pobox.spbu.ru,\\
samsonov@pink.phys.spbu.ru}
\maketitle

\begin{abstract}
The solutions of the Drinfeld equation corresponding to
the full set of different
carrier subalgebras $\frak{f}\subset\frak{sl}_3$
are explicitly constructed.
The obtained Hopf structures are studied.
It is demonstrated that the presented twist deformations can be considered
as limits of the corresponding quantum analogues ($q$-twists)
defined for the $q$-quantized algebras.
\end{abstract}

PACS numbers: 0220, 0365

AMS numbers: 17B37, 20G42


\section{Introduction}

Triangular Hopf algebras ${\cal A} (m,\Delta ,S,\eta
,\epsilon;{\cal R})$ \cite{FRT} play an essential role in quantum
theory and in particular for models with noncommutative
space-time \cite{CKNT,WESS}.
Quantizations of antisymmetric $r$-matrices (solutions of classical Yang-Baxter
equation) form an important class of such algebras.
They describe Poisson structures compatible with the initial
Lie algebra $\frak{g}$, i.e. the mechanical systems that can exist
on a space whose noncommutativity is fixed by $\frak{g}$.
Such quantum algebras can be constructed in terms of $r$-matrices by means of
Campbell-Hausdorff series \cite{D-83}. However, these
constructions are obviously inappropriate for an efficient usage
of quantum ${\cal R}$-matrices. If one provides the elements of the initial Lie
algebra $\frak{g}$ with primitive coproducts $\Delta^{\rm prim}$ and
consider the universal enveloping algebra  $U(\frak{g})$ as a Hopf algebra
with the costructure generated by $\Delta^{\rm prim}$, then the solution
${\cal F} \in U(\frak{g})^{\otimes 2}$ of the twist equation \cite{D-83}
\begin{equation}
\label{Drin-1}
{\cal F}_{12}(\Delta\otimes{\rm id})({\cal F})
={\cal F}_{23}({\rm id}\otimes\Delta)({\cal F}),
\end{equation}
\begin{equation}
\label{Drin-2}
(\epsilon\otimes{\rm id}){\cal F}=({\rm id}\otimes\epsilon){\cal F}=1
\end{equation}
allows one to find the solution of the Yang-Baxter equation, namely:
 ${\cal R}_{\cal F}= {\cal F}_{21}{\cal F}^{-1}$. Thus to obtain
the set of solutions of the twist equations (\ref{Drin-1})
(the set of twists ${\cal F}$) for the full set of carrier subalgebras
in simple Lie algebras
$\frak{g}$ used in constructing physical models is an important task
(see for example \cite{BGGS} and the references therein).
Such set of twists for a Lie algebra $\frak{g}$ will be considered
complete if for any class $\frak{h}$ of equivalent antisymmetric
solutions for the classical Yang-Baxter equation (CYBE)we can
attribute a twist
${\cal F}(\xi)$
(with a deformation parameter $\xi$) whose classical $r$-matrix,
$$
r_{\cal F}
=\frac{d}{d\xi}{\cal F (\xi)}_{21}{\cal F (\xi)}^{-1}|_{\xi=0},
$$
is a representative of the class $\frak{h}$.

Each $r$-matrix induces a dual map
that can be treated as a skew-symmetric bilinear form
$\omega _{r}$ that satisfy the condition
$$
\omega (x,[y,z])+\omega (y,[z,x])+\omega (z,[x,y])=0.
$$
Form $\omega _{r}$ is nondegenerate on the space of a subalgebra
$\frak{g}_{c}\subseteq \frak{g}$.
Such subalgebra $\frak{g}_{c} $ is called the carrier of $r$.
Subalgebras supplied with nondegenerate form $\omega _{r}$ are
called Lie quasi-Frobenius.
The classification problem for quasi-Frobenius Lie subalgebras is
far from being completed. The explicit classification is
known for some types of Lie algebras and in particular
for $\frak{sl}_3$ it was given by Stolin \cite{STO}.

The construction of the twisting elements is not only important but
also a difficult problem and
for a long time only few types of twists were known in
an explicit form \cite{R, OGI, GIA, KLM}.
In this paper we demonstrate that using
the factorization property \cite{KL,L} of twists the explicit solution
of equations (\ref{Drin-1}-\ref{Drin-2}) can be constructed for any
quasi-Frobenius subalgebra in $\frak{sl}_3$. In Section 2 the corresponding
Hopf algebras -- twist deformations $U\left(\frak{sl}_3\right)
\longrightarrow U_{F}\left(\frak{sl}_3\right)$
-- are classified and studied.

The second fundamental problem in the study of triangular Hopf algebras is the
relation between twist deformations and ordinary quantizations ($q$-deformations).
In Section 3 we study the possibility to attribute to each twist ${\cal F}$ (and
the corresponding deformation
$U_{{\cal F}}\left(\frak{f}\right)$) the quantum twist
$F_q^{\prime}$ defining deformation of quantized current algebra
$U_{q}(\frak{sl}_{3})$ and
specializing to a twist ${\cal F}$ of $U(\frak{sl}_{3})$ in the limit
$q\rightarrow 1$.
The limit is assumed to be taken along
some curve $q=1+f(s)$ as $f(s)\rightarrow 0$ for $s \rightarrow 0$.
In particular, we can obtain quantum versions of the quasi-Frobenius subalgebras
through commutativity of the following diagram
\begin{equation}
\label{q-tw}
\begin{array}{ccccc}
U\left(\frak{f}\right) & \stackrel{q}{\longrightarrow} & U_{q}\left(\frak{f}\right)
& \stackrel{\iota_{q}}{\longrightarrow}&
U_{q}(\widehat{\frak{sl}}_{3})\\
{\cal F} \downarrow &  & F_{q}\downarrow & & F^{\prime}_{q}\downarrow \\
U_{{\cal F}}\left(\frak{f}\right) &
\stackrel{q}{\longrightarrow} & U_{F_q}\left(\frak{f}\right)&
\stackrel{\iota_{q}}{\longrightarrow}&
U_{F_{q}^{\prime}}(\widehat{\frak{sl}}_{3})
\end{array}
\end{equation}
where $\iota_{q}$ is the embedding, $\widehat{\frak{sl}}_{3}$ stands for
the affinization of $\mathfrak{sl}_{3}$: $A_{2}^{(1)}$ or
$A_{2}^{(2)}$ and
$F^{\prime}_{q}:=(\iota_{q}\otimes\iota_{q})F_{q}(\iota^{-1}_{q}\otimes\iota^{-1}_{q})$.
Formally our quantization $U_{q}(\frak{sl}_{3})$ is defined over
$\mathbb{C}[s,s^{-1}]$ and the specialization $q\rightarrow~1$
corresponds to the limit $s\rightarrow 0$. In Section 3 we step by
step demonstrate that in most cases one can straightforwardly
obtain $U_{q}(\frak{f})$ as a Hopf subalgebra of $U_{q}^{{\cal
J}}(\widehat{\frak{sl}}_{3})$, the Drinfeld-Jimbo
quantization of $U(\widehat{\frak{sl}}_{3})$ with the
comultiplication deformed by some factor ${\cal J}$:
$$
\Delta_{{\cal J}}(x)={\cal J}\Delta_{q}(x){\cal J}^{-1}.
$$
$
U_{q}(\frak{f})
$ is fixed by the fact that it is a minimal Hopf subalgebra in
$U_{q}(\frak{sl}_{3})$ with the property
$$
U_{q}(\frak{f})|_{\mathbb{C}[s]}
\rightarrow U(\frak{f}),\mbox{ }s\rightarrow 0.
$$
In the particular case of $\frak{b}^{(0)}=\{H_{13}, E_{12}+E_{23}\}$
the Borel subalgebra belongs to the subalgebra
$\frak{sl}_{2}$ that is the special subalgebra of $\frak{sl}_{3}$
and we can define the corresponding $q-$quantization
of  $U(\frak{b}^{(0)})$ by embedding it into
$U_{q}(A_{2}^{(2)})$. To define the $q-$twists corresponding
to the twists described in section 2, we utilize an assumption that most
of $q-$twists can be built out of $q-$exponentials and Abelian twists.
The only twist deformation of $U(\frak{sl}_{3})$ that seems to contradict
this assumption is the one corresponding to the $r-$matrix
$
r_{\cal{JR}}=\frac 12  H_{23}\wedge E_{23}+\eta  E_{12}\wedge E_{13}
$
(see formulas (\ref{deff-jord}) in 2.6.1).
In the Appendix the special properties of the so called peripheric twists
\cite{LO} are discussed.

\section{Classification of quasi-Frobenius subalgebras in $\frak{sl}_{3}$
and twist deformations $U_{F}(\frak{sl}_{3})$}

\subsection{Abelian two-dimensional subalgebras}

We have four classes of nonequivalent
two-dimensional subalgebras \cite{STO} denoted by

$$
\begin{array}{lclcl}
\frak{h}&\Longrightarrow &H&=&\left(
{
\begin{array}{ccc}
*& 0& 0\\
0& *& 0\\
0& 0& *\\
\end{array}
}
\right)\\
\frak{h}^{(1)}&\Longrightarrow & X&=& *
\left(
{
\begin{array}{ccc}
1& 0& 0\\
0& 1& 0\\
0& 0& -2\\
\end{array}
}
\right)+
\left(
{
\begin{array}{ccc}
0& *& 0\\
0& 0& 0\\
0& 0& 0\\
\end{array}
}
\right)\\
\frak{h}^{(0,1)}&\Longrightarrow & X
&=&
\left(
{
\begin{array}{ccc}
0& *& *\\
0& 0& 0\\
0& 0& 0\\
\end{array}
}
\right)\\
\frak{h}^{(1,1)}&\Longrightarrow&
X&=&
*\left(
{
\begin{array}{ccc}
0& 1& 0\\
0& 0& 1\\
0& 0& 0\\
\end{array}
}
\right)+
\left(
{
\begin{array}{ccc}
0& 0& *\\
0& 0& 0\\
0& 0& 0\\
\end{array}
}
\right)
\end{array}
$$

\subsection{Cartan subalgebra. Case $\frak{h}$}

Any two elements $H_{1},H_{2}\in \mathbf{H}$ form a carrier subalgebra
for the so called Abelian twist (first described by Reshetikhin \cite{R}):
\begin{equation}
\label{ab-h}
\mathcal{F}_{\mathcal{R}}=\exp (\xi ^{12}H_{1}\otimes H_{2})\quad \quad \xi
^{12}\in \mathbb{C}.
\end{equation}
Let $\Lambda $ be the root system of $\frak{sl}_{3}$.
For any $\lambda \in \Lambda $
let $E_{\lambda}\in \frak{sl}_{3}$ be the element corresponding to
the root $\lambda $.
The twisting elements $\mathcal{F}_{\mathcal{R}}$ lead to the deformed
costructure:
\begin{equation}
\begin{array}{lcl}
\Delta _{\mathcal{F_{R}}}(H_{1,2}) & = & H_{1,2}\otimes 1+1\otimes H_{1,2},
\\[0.2cm]
\Delta _{\mathcal{F_{R}}}(E_{\lambda }) & = & E_{\lambda }\otimes e^{\xi
^{12}\lambda (H_{1})H_{2}}+e^{\xi ^{12}\lambda (H_{2})H_{1}}\otimes
E_{\lambda }.
\end{array}
\label{copr1}
\end{equation}


\subsection{Mixed. Case $\frak{h}^{(1)}$}

For any $E_{\lambda },\quad \lambda \in \Lambda (%
\frak{g})$ consider the Cartan element $H_{\lambda }^{\perp }$ whose
dual $\left( H_{\lambda }^{\perp }\right) ^{\ast }$
(with respect to the Killing form) is a vector orthogonal
to $\lambda $. This pair generates the two-dimensional Abelian algebra, the
carrier for the twist:
\begin{equation}
\label{ab-h1}
\mathcal{F}_{\mathcal{R}}=\exp (\xi H_{\lambda }^{\perp }\otimes E_{\lambda
}).
\end{equation}
Notice that here the parameter $\xi $ can be scaled by a similarity
transformation.
\begin{equation}  \begin{array}{lcl}
\Delta _{\mathcal{F_{R}}}(H) & = & H\otimes 1+1\otimes H - \xi \lambda
(H)H_{\lambda }^{\perp }\otimes E_{\lambda},\quad H\in \mathbf{H}\\[0.2cm]
\Delta _{\mathcal{F_{R}}}(E_{\mu })|_{\mu \neq -\lambda } & = & E_{\mu
}\otimes e^{\xi \mu (H_{\lambda }^{\perp })E_{\lambda }}+1\otimes E_{\mu
}+\xi H_{\lambda }^{\perp }\otimes E_{\lambda +\mu } \\[0.2cm]
\Delta _{\mathcal{F_{R}}}(E_{-\lambda }) & = & E_{-\lambda }\otimes
1+1\otimes E_{-\lambda }+ \\
&  & +\xi H_{\lambda }^{\perp }\otimes H_{\lambda }-\xi ^{2}(H_{\lambda
}^{\perp })^{2}\otimes E_{\lambda }
\end{array}
\label{copr2}
\end{equation}
(It is assumed that $E_{\left\{ \lambda +\mu  \right\} }$ is zero iff
$\lambda +\mu $
is not a root.)


\subsection{$\mathbf{N}^+$ subalgebra. Cases $\frak{h}^{(0,1)}$
and $\frak{h}^{(1,1)}$}

Consider, for example, the commuting generators $E_{12}+\mu E_{23},E_{13}$
and the corresponding Reshetikhin twist
\begin{equation}
\label{ab-h01-11}
\mathcal{F}_{\mathcal{R}}=\exp (\xi (E_{12}+\mu E_{23})\otimes E_{13}).
\end{equation}
Here the cases $\mu =0$ and $\mu \neq 0$ are to be considered separately.
The twist $\quad \mathcal{F}_{\mathcal{R}}\left( \xi ,\mu =0\right) \quad $%
cannot be scaled to $\quad \mathcal{F}_{\mathcal{R}}\left( \xi ,\mu \neq
0\right) $. Though it is a limit for the family \quad $\left\{ \mathcal{F}_{%
\mathcal{R}}\left( \xi ,\mu \right) \right\} \quad $as is clearly seen from
the explicit costructure below.
\begin{equation}
\begin{array}{lcl}
\nonumber
\Delta _{\mathcal{F_{R}}}(E_{12}+\mu E_{23}) & = & (E_{12}+\mu
E_{23})\otimes 1+1\otimes (E_{12}+\mu E_{23}), \rule{14mm}{0mm}
\end{array}
\end{equation}
\begin{equation}
\begin{array}{lcl}
\nonumber
\Delta _{\mathcal{F_{R}}}(E_{13}) & = & E_{13}\otimes 1+1\otimes E_{13}, \\%
[0.2cm]
\Delta _{\mathcal{F_{R}}}(H_{12}) & = & H_{12}\otimes 1+1\otimes H_{12}-3\xi
E_{12}\otimes E_{13}+\frac{3}{2}\xi ^{2}\mu E_{13}\otimes (E_{13})^{2} \\%
[0.2cm]
\Delta _{\mathcal{F_{R}}}(H_{23}) & = & H_{23}\otimes 1+1\otimes H_{23}-3\xi
\mu E_{23}\otimes E_{13}-\frac{3}{2}\xi ^{2}\mu E_{13}\otimes (E_{13})^{2} \\%
[0.2cm]
\Delta _{\mathcal{F_{R}}}(E_{21}) & = & E_{21}\otimes 1+1\otimes E_{21}+ \\%
[0.2cm]
&  & +\xi H_{12}\otimes E_{13}-\xi (E_{12}+\mu E_{23})\otimes E_{23} \\
&  & -\frac{1}{2}\xi ^{2}(2E_{12}-\mu E_{23})\otimes (E_{13})^{2}+\frac{1}{2}%
\xi ^{3}\mu E_{13}\otimes (E_{13})^{3},
\end{array}
\end{equation}
\begin{equation}
 \begin{array}{lcl}
\Delta _{\mathcal{F_{R}}}(E_{23}) & = & E_{23}\otimes 1+1\otimes E_{23}+\xi
E_{13}\otimes E_{13}, \\[0.2cm]
\Delta _{\mathcal{F_{R}}}(E_{32}) & = & E_{32}\otimes 1+1\otimes E_{32}+\xi
\mu H_{23}\otimes E_{13} \\[0.2cm]
&  & +\frac{1}{2}\xi ^{2}\mu (E_{12}-2\mu E_{23})\otimes (E_{13})^{2}+ \\%
[0.2cm]
&  & +\xi (E_{12}+\mu E_{23})\otimes E_{12}-\frac{1}{2}\xi ^{3}\mu
^{2}E_{13}\otimes (E_{13})^{3}, \\[0.2cm]
\Delta _{\mathcal{F_{R}}}(E_{31}) & = & E_{31}\otimes 1+1\otimes E_{31}+\xi
(-E_{32}+\mu E_{21})\otimes E_{13}+ \\[0.2cm]
&  & +\xi (E_{12}+\mu E_{23})\otimes H_{13}+\frac{1}{2}\xi ^{2}\mu
(H_{12}-H_{23})\otimes (E_{13})^{2} \\
&  & -\xi ^{2}(E_{12}+\mu E_{23})^{2}\otimes E_{13}+ \\[0.2cm]
&  & +\frac{1}{2}\xi ^{3}\mu (-E_{12}+\mu E_{23})\otimes (E_{13})^{3}+\frac{1%
}{2}\xi ^{4}\mu ^{2}E_{13}\otimes (E_{13})^{4},
\end{array}
\label{copr3}
\end{equation}

\subsection{$\mathbf{B}(2)$ subalgebras and Jordanian twists}

In this section the carrier algebra is normalized as
\begin{equation}
\left[ H,E\right] =E,
\end{equation}
and the twist is Jordanian \cite{OGI}:
\begin{equation}
\mathcal{F}_{\mathcal{J}}=\exp (H\otimes \sigma \left( \xi \right) ),\quad
\quad \sigma \left( \xi \right) =\ln (1+\xi E),  \label{jord}
\end{equation}
the parameter can be scaled by $\mathrm{ad}(H)$. Our task is to enumerate
the inequivalent $\mathbf{B}(2)$ subalgebras.

Choose the generator $E$ in $\mathbf{N}^{+}$. We must distinguish the cases
where the Cartan generator for $\mathbf{B}(2)$ can be diagonolized in $%
\mathbf{N}^{+}$ and where it cannot. Only two-dimensional eigenspaces (for
a Cartan subalgebra) could be found in $\mathbf{N}^{+}$. Notice that one can
always add $E$ to $H$ without principal influence on the results, we shall
not consider such cases separately.


\subsubsection{Irregular element $H$ and $E=E_{13}$}

It is sufficient to specialize $H$ as follows:
\begin{equation}
H=a^{i}E_{ii}+\alpha E_{23}+\beta E_{12},\quad \quad \sum a^{i}=0.
\end{equation}
The coefficients $\alpha $ and $\beta $ can be scaled by an appropriate $%
\mathrm{ad}(H^{\perp })$. We can have only the element proportional to $%
E=\gamma E_{13}$ here. The answer is a continuous family of $\mathbf{B}(2)$
subalgebras:
\begin{equation}
H=\mu E_{11}+(1-2\mu )E_{22}+(\mu -1)E_{33}+\alpha E_{23}+\beta E_{12},\quad
\quad E=E_{13}.
\end{equation}
Up to equivalence the corresponding Jordanian twist is represented
by the following expression
$$
{\cal F}_{irr}=\exp((H_{12}^{\perp}+\eta E_{12})\otimes
\sigma_{13}).
$$
The twisted coproducts are
$$
 \begin{array}{lll}
\Delta (H_{12}) & = & H_{12}\otimes 1-2\eta
E_{12}\otimes\sigma_{13}
-(H_{12}^{\perp}+\eta E_{12})\otimes(1-e^{-\sigma_{13}})+1\otimes
H_{12};\\
\Delta (H_{23}) & = & H_{23}\otimes 1+\eta
E_{12}\otimes\sigma_{13}
-(H_{12}^{\perp}+\eta E_{12})\otimes(1-e^{-\sigma_{13}})+1\otimes
H_{23};\\
\Delta (E_{12}) & = & E_{12}\otimes 1+1\otimes
E_{12};\\
\Delta (E_{23}) & = & E_{23}\otimes e^{\sigma _{13}}+ \eta
E_{13}\otimes \sigma_{13} e^{\sigma_{13}}+1\otimes
E_{23};\\
\Delta (E_{13}) & = & E_{13}\otimes e^{\sigma _{13}}+1\otimes
E_{13};\\
\Delta(E_{32})&=& E_{32}\otimes e^{-\sigma_{13}}+
(H_{12}^{\perp}+\eta E_{12})\otimes
E_{12}e^{-\sigma_{13}}+1\otimes E_{32};\\
 \Delta(E_{21})&=&E_{21}\otimes 1+\eta H_{12}\otimes\sigma_{13}
 -\eta^{2}E_{12}\otimes(\sigma_{13})^{2}\\
 &&-(H_{12}^{\perp}+\eta E_{12})\otimes E_{23}e^{-\sigma_{13}}
 +1\otimes E_{21};\\
 \Delta(E_{31})&=& E_{31}\otimes e^{-\sigma_{31}}-\eta
 E_{32}\otimes\sigma_{13}e^{-\sigma_{13}}+
 (H_{12}^{\perp}+\eta E_{12})\otimes H_{13}e^{-\sigma_{13}}+\\
 &&+((H_{12}^{\perp}+\eta E_{12})
 -(H_{12}^{\perp}+\eta E_{12})^{2})\otimes
 (e^{-\sigma_{13}}-e^{-2\sigma_{13}})+1\otimes E_{31};
\end{array}
$$
Thus the irregularity results in the appearance of
$\sigma_{13}$ in the coproducts.
When both parameters $\alpha $ and $\beta $\ are zeros we come to the
regular case treated below.


\subsubsection{Regular $H$}

Let
\begin{equation}
H=a^{i}E_{ii},\quad \quad \sum a^{i}=0.
\end{equation}
In the case $E=E_{13}$ all the properties of the Jordanian twist are the
same as described above.

\begin{enumerate}
\item  $E=E_{\lambda }$, with $\lambda $ being one of the simple roots.
Let $E=E_{12}$. Again we find the parameterized family of subalgebras (and
correspondingly the twists):
\begin{equation}
\left\{ H=\mu E_{11}+(\mu -1)E_{22}+(1-2\mu )E_{33},\quad E_{12}\right\} .
\label{reg-1}
\end{equation}
The case $E=E_{23}$ is treated analogously.

\item  $E=E_{12}+E_{23}$ . The case $E=E_{12}+\gamma E_{23}$ can be scaled
to the normalized one: $E=E_{12}+E_{23}$. Immediately we find $%
a^{1}=1,a^{2}=0$, thus:
\begin{equation}
\left\{ H=E_{11}-E_{33},\quad \quad E=E_{12}+E_{23}\right\} .
\label{reg-2}
\end{equation}

\item  $E=E_{12}+E_{13}$. Again only the normalized combination is to be
considered. The algebra is unique:
\begin{equation}
\left\{ H=\frac{2}{3}E_{11}-\frac{1}{3}E_{22}-\frac{1}{3}E_{33},\quad \quad
E=E_{12}+E_{13}\right\} .
\end{equation}
The alternative combination $E=E_{23}+E_{13}$ is treated analogously.
\end{enumerate}


\subsection{4-dimensional carriers}

\subsubsection{One Cartan generator}

In this case we can assume that the carrier belongs to the Borel subalgebra:
$\mathbf{L}\subset \mathbf{B}$.  Let  $H_{\lambda =e_{1}-e_{3}}^{\perp }\equiv
H_{13}^{\perp }$ (with the dual vector $(H_{13}^{\perp })^{\ast }$
orthogonal to the highest root $e_{1}-e_{3}$) and assume that $a^{1}\neq
a^{3}$ in $H=a^{i}E_{ii}$.

In this case the family of carrier algebras $L\left( \alpha ,\beta ,\gamma
,\delta \right) \subset \mathbf{B}$\ ,
\begin{equation}
\begin{array}{ll}
\left[ H,A\right] =\alpha A, & \left[ H,E\right] =\delta E, \\
\left[ H,B\right] =\beta B, & \left[ A,B\right] =\gamma E, \\
\alpha +\beta =\delta . & \alpha ,\beta ,\gamma ,\delta \in \mathbb{C}
\end{array}
\label{l-carrier}
\end{equation}
a representative can be chosen with $\delta =1$ while $\gamma $ will finally
coincide with the deformation parameter. So we are to consider only the case
$L\left( \alpha ,\beta \right) $:
\begin{equation}
\begin{array}{ll}
\left[ H,A\right] =\alpha A, & \left[ H,E\right] =E, \\
\left[ H,B\right] =\beta B, & \left[ A,B\right] =E, \\
\alpha +\beta =1. &
\end{array}
\label{l-rel}
\end{equation}
The set $\left\{ \mathbf{L}_{\alpha ,\beta }\right\} $ of carriers (\ref
{l-rel}) is to be further classified due to the values of the second
cohomology group $H^{2}(\mathbf{L},\mathbb{C})$ and the orbits of the
normalizer $N(\mathbf{L})$ (of $\mathbf{L}$ in $\frak{sl}_{3}$)
 formed by its adjoint action in $H^{2}(\mathbf{L},\mathbb{C})$ \cite{STO}.
 Consider the list of cohomological
properties of $\mathbf{L}$:
\begin{equation}
\begin{array}{cccc}
& \mathrm{dim}Z^{2} & \mathrm{dim}B^{2} & \mathrm{dim}H^{2} \\
\alpha ,\beta \neq 0,-1 & 3 & 3 & 0 \\
\left.
\begin{array}{c}
\alpha =0,\beta =1 \\
\alpha =1,\beta =0
\end{array}
\right\}  & 3 & 2 & 1 \\
\left.
\begin{array}{c}
\alpha =-1,\beta =2 \\
\alpha =2,\beta =-1
\end{array}
\right\}  & 4 & 3 & 1
\end{array}
\label{cohomol}
\end{equation}
The first column describes the subsets in $\left\{ \mathbf{L}_{\alpha ,\beta
}\right\} $ that are to be considered separately.

\begin{enumerate}
\item  The case $\alpha ,\beta \neq 0,-1$. The carrier is the Frobenius
subalgebra with the nondegenerate coboundary $\omega =E^{\ast }([,])$. The
corresponding twist is the extended Jordanian twist \cite{KLM}, it can be
written in two (equivalent) forms:
\begin{equation}
\begin{array}{l}
\mathcal{F}_{\mathcal{E}}=\exp (\xi A\otimes Be^{-\beta \sigma (\xi )})\exp
(H\otimes \sigma (\xi )), \\
\mathcal{F}_{\mathcal{E}^{\prime }}=\exp (-\xi B\otimes Ae^{-\alpha \sigma
(\xi )})\exp (H\otimes \sigma (\xi )).
\end{array}
\label{norm-ext}
\end{equation}

They are connected by the automorphism
\[
\begin{array}{lll}
i & : & \left\{
\begin{array}{lll}
A & \longrightarrow  & -B \\
B & \longrightarrow  & A \\
\alpha  & \rightleftarrows  & \beta
\end{array}
\right\} .
\end{array}
\]

For example in the case of $\mathcal{F}_{\mathcal{E'}}$ the costructure is
defined by the relations
\begin{equation}
\begin{array}{lcl}
\Delta _{\mathcal{E}}(H) & = & H\otimes e^{-\sigma (\xi )}+1\otimes H-\xi
A\otimes Be^{-(\beta +1)\sigma (\xi )}, \\[0.2cm]
\Delta _{\mathcal{E}}(A) & = & A\otimes e^{-\beta \sigma (\xi )}+1\otimes A,
\\[0.2cm]
\Delta _{\mathcal{E}}(B) & = & B\otimes e^{\beta \sigma (\xi )}+e^{\sigma
(\xi )}\otimes B, \\[0.2cm]
\Delta _{\mathcal{E}}(E) & = & E\otimes e^{\sigma (\xi )}+1\otimes E.
\end{array}
\label{delf}
\end{equation}
The $\mathcal{R}$-matrix has the form
\begin{equation}
\begin{array}{lcl}
\mathcal{R}_{\mathcal{E}} & = & \exp (\xi Be^{-\beta \sigma (\xi )}\otimes
A)\exp (\sigma (\xi )\otimes H) \\
&  & \exp (-H\otimes \sigma (\xi ))\exp (-\xi A\otimes Be^{-\beta \sigma
(\xi )}) \\
& = & 1\otimes 1-\xi r_{\mathcal{E}}+\mathcal{O}(\xi ^{2}).
\end{array}
\end{equation}
The corresponding classical $r$-matrix is
\begin{equation}
r_{\mathcal{E}}=H\wedge E+A\wedge B.
\end{equation}
The parameter $\xi $ in (\ref{norm-ext}) can be scaled. Notice that equal
coefficients in two terms of $r_{\mathcal{E}}$-matrix is
the necessary and sufficient condition for the corresponding
form $\omega _{\mathcal{E}}$ to be a cocycle.
The reason is that the extension factors in $%
\mathcal{F}_{\mathcal{E}}$ and $\mathcal{F}_{\mathcal{E}^{\prime }}$ are the
discrete twists and can only borrow the continuous parameter from the smooth
set of Hopf algebras (twisted by $\mathcal{F}_{\mathcal{J}}=\exp (H
\otimes \sigma (\xi ))$).

When the carrier $\mathbf{L}_{\alpha ,\beta }$ is identified with the
subalgebra of $sl(3)$ it is convenient to describe the freedom in its
definition by introducing the second Cartan generator $H_{13}^{\perp }=\frac{%
1}{3}E_{11}-\frac{2}{3}E_{22}+\frac{1}{3}E_{33}$ and the parameter $\zeta
=\alpha -\frac{1}{2}$. In these terms the twisting element
$\mathcal{F}_{\mathcal{E}}$ from (\ref{norm-ext})
takes the form
\[
\mathcal{F}_{\mathcal{E}}=\exp \left(\xi E_{12}\otimes E_{23}e^{\left( \zeta -%
\frac{1}{2}\right) \sigma \left(\xi \right)}\right)\exp \left(\left(\frac 12
H_{13}+\zeta H_{13}^{\perp }\right) \otimes \sigma (\xi )\right).
\]
The deformed $U_{\mathcal{E}}\left( sl(3\right) $ is defined by the
following coproducts:
\begin{equation}
\begin{array}{lll}
 \Delta_{\mathcal{E}}(H_{12})&=& H_{12}\otimes 1 +1\otimes H_{12}
+(\frac 12 H_{13}^{\perp}+\zeta H_{13}^{\perp})\otimes
(e^{-\sigma\left(\xi\right)}-1)-\\
&&-\xi E_{12}\otimes E_{23}
e^{(\zeta-\frac 32)\sigma(\xi)};\\
\Delta _{\mathcal{E}}(H_{23}) & = & H_{23}\otimes 1+1\otimes
H_{23}+\left(\frac 12 H_{13}+\zeta H_{13}^{\perp }\right) \otimes
\left(e^{-\sigma
\left( \xi \right)}-1\right)-\\
&  & -\xi E_{12}\otimes E_{23}e^{\left(\zeta
-\frac{3}{2}\right)\sigma\left( \xi \right)}; \\
\Delta _{\mathcal{E}}(E_{12}) & = & E_{12}\otimes e^{^{\left( \zeta
-\frac{1%
}{2}\right) \sigma\left( \xi \right) }}+1\otimes E_{12}; \\
\Delta _{\mathcal{E}}(E_{23}) & = & E_{23}\otimes e^{\left( \frac{1}{2}%
-\zeta \right) \sigma\left( \xi \right) }+e^{\sigma\left( \xi
\right) }\otimes E_{23}; \\
\Delta _{\mathcal{E}}(E_{13}) & = & E_{13}\otimes e^{\sigma \left(
\xi
\right) }+1\otimes E_{13}; \\
\Delta _{\mathcal{E}}(E_{21}) & = & E_{21}\otimes e^{-\left( \frac
12+\zeta \right)
\sigma\left( \xi \right) }+1\otimes E_{21}+ \\
&  & +\xi \left( H_{12}-\frac 12 H_{13}-\zeta H_{13}^{\bot
}\right) \otimes
E_{23}e^{-\sigma\left( \xi \right) }; \\
\Delta _{\mathcal{E}}(E_{32}) & = & E_{32}\otimes e^{^{\left( \zeta
-\frac{1%
}{2}\right) \sigma\left( \xi \right) }}+1\otimes E_{32}
+\xi E_{12}\otimes H_{23}e^{^{\left( \zeta
-\frac{1}{2}\right) \sigma\left( \xi \right) }}+ \\
&  & +\xi
\left(\frac 12 H_{13}+\zeta H_{13}^{\bot }\right) \otimes
E_{12}e^{-\sigma
\left( \xi \right) }-\\
&  &-\xi\left(\frac 12 H_{13}+\zeta H_{13}^{\perp }\right)
E_{12}\otimes\left( e^{\left(\zeta-\frac 12\right)\sigma\left( \xi
\right)}-e^{\left(\zeta-\frac 32\right)\sigma \left(
\xi\right)}\right)\\
&  & -\xi^{2}E_{12}\otimes E_{23}E_{12}e^{^{\left( \zeta -\frac{3}{2}%
\right) \sigma\left( \xi \right) }}-\xi^{2} E_{12}^{2}\otimes
E_{23}e^{2\left(\zeta-\frac 12\right)\sigma
\left(\xi\right)}; \\
\end{array}
\end{equation}
\begin{equation}
\begin{array}{ll}
\Delta _{\mathcal{E}}(E_{31}) & =  E_{31}\otimes e^{-\sigma
\left( \xi
\right) }+1\otimes E_{31}+ \\
& + \xi \left(\frac 12 H_{13}+\zeta H_{13}^{\perp }\right)
\otimes
H_{13}e^{-\sigma\left( \xi \right) }+\\
& + \xi \left( 1-\frac 12 H_{13}-\zeta H_{13}^{\perp }\right)
\left(\frac 12 H_{13}+\zeta H_{13}^{\perp }\right) \otimes \left(
e^{-\sigma\left( \xi \right)}-e^{-2\sigma\left( \xi \right) }\right) + \\
&  + \xi ^{2}\left(\frac 12 H_{13}+\zeta H_{13}^{\perp }-1\right)
E_{12}\otimes E_{23}\left( e^{^{\left( \zeta -\frac{3}{2}\right)
\sigma\left( \xi \right) }}-2e^{^{\left( \zeta -\frac{5}{2}\right)
\sigma\left( \xi
\right) }}\right) \\
&   +\xi E_{12}\otimes E_{21}e^{^{\left( \zeta
-\frac{1}{2}\right) \sigma
\left( \xi \right) }}-\xi E_{32}\otimes E_{23}e^{^{\left( \zeta -\frac{3%
}{2}\right) \sigma\left( \xi \right) }}- \\
&   -\xi ^{2}E_{12}\otimes H_{13}E_{23}e^{^{\left( \zeta -\frac{3}{2}%
\right) \sigma\left( \xi \right) }}+\xi ^{3}E_{12}^{2}\otimes
E_{23}^{2}e^{^{2\left( \zeta -\frac{3}{2}\right) \sigma\left( \xi
\right) }};
\end{array}
\label{gen-ext-jord}
\end{equation}

\item  The case $\alpha =-1,\beta =2$. There are two twists for the carrier
algebra $\mathbf{L}(-1,2)$.

First we have the coboundary form of the previous type $\omega =E^{\ast
}([,])$. And the corresponding twists (\ref{norm-ext}) with $\alpha
=-1,\beta =2$.

The second possibility is due to the nontrivial elements of the cohomology
group $H^{2}\left( \mathbf{L}(-1,2)\right) $. The cochain
\[
(\psi :\mathbf{L}\wedge \mathbf{L}\longrightarrow \mathbb{C})
\]
such that
\[
\psi (A,E)\neq 0
\]
is not cohomologous to zero. This means that the form
\begin{equation}
\omega =B^{\ast }([,])+\zeta A^{\ast }\wedge E^{\ast }
\end{equation}
is a nontrivial cocycle for any $\zeta \in \mathbb{C}$. The corresponding twist
is a composition of a Reshetikhin and deformed Jordanian \cite{KL} factors.
\begin{equation}
\begin{array}{l}
\mathcal{F}_{\mathcal{JR}}=\exp (\frac{1}{2}H\otimes {\sigma (\zeta ,\xi )}%
)\exp (\zeta A\otimes E), \\
{\sigma (\zeta ,\xi )}=\ln (1+\xi B-\frac{1}{2}\xi \zeta ^{2}E^{2}).
\end{array}
\label{deff-jord}
\end{equation}
Notice that the composition is possible due to the fact that after applying
the Reshetikhin factor $\mathcal{F}_{\mathcal{R}}=\exp (\zeta A\otimes E)$
we get the primitive Borel subalgebra on the space generated by $\left\{ H,B-%
\frac{1}{2}\zeta ^{2}E^{2}\right\} $.

Here we have the universal $\mathcal{R}$-matrix:
\begin{equation}
  \mathcal{R}_{\mathcal{JR}}=\exp (\frac{1}{2}{\sigma (\zeta ,\xi )}\otimes
H)\exp (\xi E\otimes A)\exp (-\xi A\otimes E)\exp (-\frac{1}{2}H\otimes {%
\sigma (\zeta ,\xi )}).
\end{equation}
Choosing $\zeta =\eta \xi $ we get the $r$-matrix
\begin{equation}
r_{\mathcal{JR}}=\frac{1}{2}H\wedge B+\eta A\wedge E.  \label{r-jr}
\end{equation}
Obviously the term $\eta A^{\ast }\wedge E^{\ast }$ can give the
nondegenerate cocycle also with the second basic coboundary $E^{\ast }([,])$%
. The corresponding $r$-matrix has the form
\[
r_{\mathcal{JR}}^{\prime }=H\wedge E+A\wedge B-\eta H\wedge B.
\]
Nevertheless a simple substitution ($H\rightarrow H+A,\quad A\rightarrow
\eta A,\quad B\rightarrow -(1/\eta )(B+E),\quad E\rightarrow -E$) brings us
again to the $r$-matrix (\ref{r-jr}). Thus we have only two different
solutions here. For the first of them the deformed costructure can be easily
obtained as a special case of (\ref{gen-ext-jord}). To present the necessary
coproducts for the second case let us use the following injection in $sl(3)$%
:
$$
\begin{array}{l}
H=H_{23};\quad E=E_{13};\quad A=E_{12};\quad B=E_{23}; \\
{\sigma (\zeta ,\xi )}=\ln (1+\xi E_{23}-\frac{1}{2}\xi \zeta ^{2}E_{13}^{2})
\end{array}
$$
In these terms the coproducts $\Delta _{\mathcal{JR}}$ are defined by the
formulae
\begin{equation}
 \begin{array}{lll}
\Delta _{\mathcal{JR}}(H_{12}) & = & H_{12}\otimes 1+\frac{1}{2}\xi
H_{23}\otimes \left( E_{23}+\zeta E_{13}^{2}\right) e^{-{\sigma (\xi ,\zeta )%
}}- \\
&  & -3\zeta E_{12}\otimes E_{13}e^{-\frac{1}{2}{\sigma (\xi ,\zeta )}%
}+1\otimes H_{12}; \\
\Delta _{\mathcal{JR}}(H_{23}) & = & H_{23}\otimes e^{-{\sigma (\xi ,\zeta )}%
}+1\otimes H_{23}; \\
\Delta _{\mathcal{JR}}(E_{12}) & = & E_{12}\otimes e^{-\frac{1}{2}{\sigma
(\xi ,\zeta )}}-\frac{1}{2}\xi H_{23}\otimes E_{13}e^{-{\sigma (\xi ,\zeta )}%
}+1\otimes E_{12}; \\
\Delta _{\mathcal{JR}}(E_{23}) & = & E_{23}\otimes e^{{\sigma (\xi ,\zeta )}%
}+\zeta E_{13}\otimes E_{13}e^{\frac{1}{2}{\sigma (\xi ,\zeta )}}+1\otimes
E_{23}; \\
\Delta _{\mathcal{JR}}(E_{13}) & = & E_{13}\otimes e^{\frac{1}{2}{\sigma
(\xi ,\zeta )}}+1\otimes E_{13}; \\
\Delta _{\mathcal{JR}}(E_{21}) & = & E_{21}\otimes e^{\frac{1}{2}{\sigma
(\xi ,\zeta )}}-\zeta E_{12}\otimes \left( E_{23}+\zeta E_{13}^{2}\right)
e^{-\frac{1}{2}{\sigma (\xi ,\zeta )}}+ \\
&  & +\zeta H_{12}\otimes E_{13}+\frac{1}{2}\xi \zeta H_{23}\otimes
E_{13}E_{23}e^{-{\sigma (\xi ,\zeta )}}+1\otimes E_{21}; \\
\Delta _{\mathcal{JR}}(E_{32}) & = & E_{32}\otimes e^{-\sigma (\xi ,\zeta
)}+1\otimes E_{32}+ \\
&  & +\zeta E_{12}\otimes E_{12}e^{-\frac{1}{2}\sigma (\xi ,\zeta )}-\frac{1%
}{2}\xi \zeta H_{23}E_{12}\otimes E_{13}e^{-\frac{3}{2}\sigma (\xi ,\zeta )}-
\\
&  & +\frac{1}{2}\xi H_{23}\otimes (H_{23}-\zeta E_{12}E_{13})e^{-\sigma
(\xi ,\zeta )} \\
&  & +\xi ^{2}(\frac{1}{2}H_{23}-\frac{1}{4}H_{23}^{2})\otimes (E_{23}-\zeta
E_{13}^{2})e^{-2\sigma (\xi ,\zeta )}; \\
\Delta _{\mathcal{JR}}(E_{31}) & = & E_{31}\otimes e^{-\frac{1}{2}\sigma
(\xi ,\zeta )}+1\otimes E_{31}-\zeta E_{32}\otimes E_{13}e^{-\sigma (\xi
,\zeta )}+ \\
&  & +\zeta E_{12}\otimes H_{13}e^{-\frac{1}{2}\sigma (\xi ,\zeta )}-\zeta
^{2}E_{12}^{2}\otimes E_{13}e^{-\sigma (\xi ,\zeta )}+ \\
&  & +\frac{1}{2}\xi H_{23}\otimes (E_{21}+\zeta E_{13}-\zeta
H_{13}E_{13})e^{-\sigma (\xi ,\zeta )}+ \\
&  & +\frac{1}{2}\xi \zeta H_{23}E_{12}\otimes (-E_{23}+2\zeta
E_{13}^{2})e^{-\frac{3}{2}\sigma (\xi ,\zeta )}+ \\
&  & -\xi ^{2}\zeta ^{2}(\frac{1}{2}H_{23}-\frac{1}{4}H_{23}^{2})\otimes
(E_{23}E_{13}-\zeta E_{13}^{3})e^{-2\sigma (\xi ,\zeta )};
\end{array}
\end{equation}

\item  The case $\alpha =0,\beta =1$. This is the so called peripheric case
\cite{LO}. The carrier algebra $\mathbf{L}_{(0,1)}$ is defined by the relations
(\ref{l-carrier}) with $\alpha =0,\beta =1$. Again we can use the same
coboundary form $E^{\ast }([,])$ as in the case 1 and get the peripheric
versions of the twists (\ref{norm-ext}):
\begin{equation}
\begin{array}{l}
\mathcal{F}_{\mathcal{P}}=\exp (\xi A\otimes Be^{-\sigma (\xi )})\exp
(H\otimes \sigma (\xi )), \\
\mathcal{F}_{\mathcal{P}^{\prime }}=\exp (-\xi B\otimes A)\exp (H\otimes
\sigma (\xi ))
\end{array}
\label{per-ext}
\end{equation}
with the costructure (for the version $\mathcal{F}_{\mathcal{P}}$)
\begin{equation}
\begin{array}{lcl}
\Delta _{\mathcal{P}}(H) & = & H\otimes e^{-\sigma (\xi )}+1\otimes H-\xi
A\otimes Be^{-2\sigma (\xi )}, \\[0.2cm]
\Delta _{\mathcal{P}}(A) & = & A\otimes e^{\sigma (\xi )}+1\otimes A, \\%
[0.2cm]
\Delta _{\mathcal{P}}(B) & = & B\otimes e^{\sigma (\xi )}+e^{\sigma (\xi
)}\otimes B, \\[0.2cm]
\Delta _{\mathcal{P}}(E) & = & E\otimes e^{\sigma (\xi )}+1\otimes E.
\end{array}
\label{delp}
\end{equation}
We have also the cohomologically nontrivial map $\omega $ that can be chosen
to be
\begin{equation}
\omega _{H}=H^{\ast }\wedge A^{\ast }.
\end{equation}
The only coboundary map that can extend this $\omega _{H}$ to create a
nondegenerate form for $\mathbf{L}(0,1)$ is again $E^{\ast }([,])$,
\begin{equation}
\omega =E^{\ast }([,])+\zeta \omega _{H}=E^{\ast }([,])+\zeta H^{\ast
}\wedge A^{\ast }.  \label{om-rp}
\end{equation}
The inverse of the $\omega $-form matrix acquires the additional term
proportional to $B\wedge E$. Notice that the costructure (\ref{delp})
provides a pair of commuting primitive elements: $\sigma $ and $Be^{-\sigma }
$. This signifies the possibility to apply the corresponding Reshetikhin
twist to the algebra $U_{\mathcal{P}}(\mathbf{L}(0,1))$ deformed by (\ref
{per-ext}). Again for the version $\mathcal{F}_{\mathcal{P}}$ we have the
composition:
\begin{equation}
\mathcal{F}_{\mathcal{RP}}=\exp (\eta Be^{-\sigma (\xi )}\otimes \sigma (\xi
))\exp (\xi A\otimes Be^{-\sigma (\xi )})\exp (H\otimes \sigma (\xi )).
\label{per-resh}
\end{equation}
Parameters $\xi $ and $\eta $ are independent. Putting $\eta =\zeta \xi $ we
arrive at the $r$-matrix
\begin{equation}
r_{\mathcal{RP}}=H\wedge E+A\wedge B+\zeta B\wedge E
\end{equation}
which is in accord with $\omega $ in (\ref{om-rp}). This construction can be
easily implemented for the case $\alpha =1,\beta =0$ with similar results.
To conclude this point we must add that the two basic coboundary maps $%
E^{\ast }([,])$ and $B^{\ast }([,])$ can certainly be combined. This
corresponds to the redefinition of the extension in the basic peripheric
twists:
\begin{equation}
\begin{array}{l}
\mathcal{F}_{\mathcal{P}}=\exp (\xi A\otimes (B+\zeta E)e^{-\sigma (\xi
)})\exp (H\otimes \sigma (\xi )), \\
\mathcal{F}_{\mathcal{P}}=\exp (-\xi (B+\zeta E)\otimes A)\exp (H\otimes
\sigma (\xi )).
\end{array}
\label{per-ext2}
\end{equation}
The latter is possible due to the equal eigenvalues of $\mathrm{ad}(H)$ on $B
$ and $E$.
\end{enumerate}


\subsubsection{Two Cartan generators}

For any carrier of the type $\mathbf{L}_{\alpha, \beta}$ one can find in $%
\frak{g}$ the element $H^{\perp}$ (that commutes with $E$) and as a
result remains primitive after the twist $\mathcal{F}_{\mathcal{E}}$ or $%
\mathcal{F}_{\mathcal{P}}$. Notice that in this case the extended twist
multiplied by the additional Reshetikhin factor $e^{\zeta H^{\perp} \otimes
\sigma}$ is equivalent to the shift of the Cartan element, $H
\longrightarrow H + \zeta H^{\perp}$, in the initial extended twist.
(When the additional twisting element is dragged  through the extension factor
 the power in the exponent $e^{- \beta \sigma}$ is changed  because $\beta$ is
 shifted together with $H$.) Thus the additional factor does not produce new
twist but results in changes of parameters of the carrier $\mathbf{L}%
_{\alpha, \beta}$.

On the contrary when two commuting elements are taken from $\mathbf{N}^{+}$
(for example $E_{12}$ and $E_{13}$) the Cartan elements can be chosen so
that the four-dimensional carrier algebra will have the structure of a
direct sum of two $\mathbf{B}(2)$ subalgebras:
\begin{equation}
\left[ H_{12}^{\perp },E_{13}\right] =E_{13},\quad \quad \left[
H_{13}^{\perp },E_{12}\right] =E_{12}.  \label{two-bor}
\end{equation}
with
\[
\begin{array}{ccl}
H_{12}^{\perp } & = & \frac{1}{3}(E_{11}+E_{22})-\frac{2}{3}E_{33}, \\[1ex]
H_{13}^{\perp } & = & \frac{1}{3}(E_{11}+E_{33})-\frac{2}{3}E_{22}.
\end{array}
\]
Both $\mathbf{B}(2)$ subalgebras can be twisted by Jordanian twists
simultaneously with independent parameters, both can be scaled (just as in
the case of unique Jordanian twist).
\begin{equation}
\mathcal{F}_{\mathcal{JJ}}=\exp (H_{13}^{\perp }\otimes \sigma _{12}(\xi
_{1}))\exp (H_{12}^{\perp }\otimes \sigma _{13}(\xi _{2})).  \label{two-cart}
\end{equation}
\begin{equation*}
\begin{array}{ccl}
\Delta _{\mathcal{JJ}}(H_{12}^{\perp }) & = & H_{12}^{\perp }\otimes
e^{-\sigma _{13}(\xi _{2})}+1\otimes H_{12}^{\perp }; \\
\Delta _{\mathcal{JJ}}(H_{13}^{\perp }) & = & H_{13}^{\perp }\otimes
e^{-\sigma _{12}(\xi _{1})}+1\otimes H_{13}^{\perp }; \\
\Delta _{\mathcal{JJ}}(E_{12}) & = & E_{12}\otimes e^{\sigma _{12}(\xi
_{1})}+1\otimes E_{12}; \\
\Delta _{\mathcal{JJ}}(E_{13}) & = & E_{13}\otimes e^{\sigma _{13}(\xi
_{2})}+1\otimes E_{13}; \\
\Delta _{\mathcal{JJ}}(E_{23}) & = & E_{23}\otimes e^{-\sigma _{12}(\xi
_{1})}e^{\sigma _{13}(\xi _{2})}+\xi H_{13}^{\perp }\otimes E_{13}e^{-\sigma
_{12}(\xi _{1})}+1\otimes E_{23}; \\
\Delta _{\mathcal{JJ}}(E_{21}) & = & E_{21}\otimes e^{-\sigma _{12}(\xi
_{1})}+1\otimes E_{21}- \\
&  & -\xi _{2}H_{12}^{\perp }\otimes E_{23}e^{-\sigma _{13}(\xi _{2})}+\xi
_{1}H_{13}^{\perp }\otimes H_{12}e^{-\sigma _{12}(\xi _{1})}- \\
&  & -\xi _{1}H_{12}^{\perp }H_{13}^{\perp }\otimes (e^{-\sigma _{12}(\xi
_{1})}-e^{-\sigma _{12}(\xi _{1})}e^{-\sigma _{13}(\xi _{2})})+ \\
&  & +\xi _{1}(H_{13}^{\perp }-(H_{13}^{\perp })^{2})\otimes (e^{-\sigma
_{12}(\xi _{1})}-e^{-2\sigma _{12}(\xi _{1})}); \\
\Delta _{\mathcal{JJ}}(E_{32}) & = & E_{32}\otimes e^{\sigma _{12}(\xi
_{1})}e^{-\sigma _{13}(\xi _{2})}+\xi _{2}H_{12}^{\perp }\otimes
E_{12}e^{-\sigma _{13}(\xi _{2})}+1\otimes E_{32}; \\
\Delta _{\mathcal{JJ}}(E_{31}) & = & E_{31}\otimes e^{-\sigma _{13}(\xi
_{2})}+1\otimes E_{31}+ \\
&  & +\xi _{2}H_{12}^{\perp }\otimes H_{13}e^{-\sigma _{13}(\xi _{2})}-\xi
_{1}H_{13}^{\perp }\otimes E_{32}e^{-\sigma _{12}(\xi _{1})}+ \\
&  & +\xi _{2}(H_{12}^{\perp }-(H_{12}^{\perp })^{2})\otimes (e^{-\sigma
_{13}(\xi _{2})}-e^{-2\sigma _{13}(\xi _{2})})+ \\
&  & +\xi _{2}H_{12}^{\perp }H_{13}^{\perp }\otimes (e^{-\sigma _{12}(\xi
_{1})}e^{-\sigma _{13}(\xi _{2})}-e^{-\sigma _{13}(\xi _{2})});
\end{array}
\end{equation*}
Due to the appearance of two primitive $\sigma $'s (see the second two lines
of the list) the general form in this case must contain additional
Reshetikhin twist:
\begin{equation}
\mathcal{F}_{\mathcal{JJ}}=\exp (\zeta \sigma _{12}(\xi _{1})\otimes \sigma
_{13}(\xi _{2}))\exp (H_{13}^{\perp }\otimes \sigma _{12}(\xi _{1}))\exp
(H_{12}^{\perp }\otimes \sigma _{13}(\xi _{2})).  \label{two-cart-r}
\end{equation}
Only two of three parameters can be scaled here (to get a nontrivial
contribution to the $r$-matrix the parameter $\zeta $ can be chosen
proportional to $1/\xi $ with $\xi _{1}=\alpha _{1}\xi ,\xi _{2}=\alpha
_{2}\xi $).


\subsection{6-dimensional carrier}

Up to the renumeration of the basic elements there is only one
six-dimensional Frobenius subalgebra $\frak{P}$ in $\frak{sl}(3)$ with
the generators:
\begin{equation}
\left\{ H_{13}^{\perp },H_{23}^{\perp },E_{12},E_{13},E_{23},E_{32}\right\}
\end{equation}
This subalgebra can be considered as the simplest case of parabolic
subalgebras in the classical series $\frak{sl}(N)$. The parabolic subalgebra arise
when some negative simple root generators are dropped from the Chevalle basis
of a simple Lie algebra. In our case this happens when the generator $%
E_{(e_{2}-e_{1})}=E_{21}$ is eliminated from the basis. The remaining
elements generate $\frak{P}$. It is easy to check that this algebra has
a trivial second cohomology, $H^{2}(\frak{P},\mathbb{C})=0$. Thus we have one
solution, it is called the parabolic twist \cite{LS}.
\begin{equation}
\begin{array}{lcl}
\mathcal{F}_{\wp } & = &
\mathcal{F}_{\mathcal{D}}\mathcal{F}_{\mathcal{EJ}}
\\
& = & \exp (H_{13}^{\perp}\otimes (2\sigma _{13}(\xi_{1})+\sigma
_{32}(\xi_{2}))\cdot\\
& & \cdot\exp (-\xi_{1}E_{23}\otimes E_{12}e^{\sigma
_{13}(\xi_{1})})\exp (H_{23}\otimes \sigma
_{13}(\xi_{1}))\\[2ex]
&=& \exp(H_{13}^{\perp}\otimes\sigma_{13}(\xi_{1}))
\exp(H_{13}^{\perp}\otimes\sigma_{32}(\xi_{2}))
 \exp(H_{13}^{\perp}\otimes\sigma_{13}(\xi_{1}))\cdot
 {\cal F}_{\cal EJ}
\end{array}
\label{par-twist}
\end{equation}
The parabolic twist factorizes into the ordinary extended Jordanian $\mathcal{F}_{%
\mathcal{EJ}}$ and the factor $\mathcal{F}_{\mathcal{D}}$. The latter can be
considered as a deformed version of the Jordanian twist. The final deformed
costructure in $\frak{P}$ looks as follows:
\begin{equation}
 \begin{array}{ll}
\nonumber
\Delta _{\wp }(H_{13}^{\perp }) = &  (H_{13}^{\perp }\otimes
1)(1\otimes 1+\xi _{2}C(\xi _{1}))^{-1}+1\otimes H_{13}^{\perp };
\\[1ex]
\Delta _{\wp }(H_{23}^{\perp }) = &  H_{23}^{\perp }\otimes e^{-\sigma
_{13}(\xi _{1})}+1\otimes H_{232}^{\perp }+ \\
&   +\xi _{1}(E_{23}+\xi H_{12}^{\perp }H_{13}^{\perp })\otimes
E_{12}e^{-\sigma _{13}(\xi _{1})}e^{-\sigma _{32}(\xi _{2})}e^{-\sigma
_{13}(\xi _{1})};
\\[1ex]
\Delta _{\wp }(E_{12}) = &  E_{12}\otimes e^{\sigma _{32}(\xi
_{2})}e^{\sigma _{13}(\xi _{1})}+e^{\sigma _{13}(\xi _{1})}\otimes
E_{12}+\xi _{1}\xi _{2}H_{13}^{\perp }E_{12}\otimes E_{12}; \\[1ex]
\Delta _{\wp }(E_{{13}}) = &  \left(
\begin{array}{l}
E_{13}\otimes e^{\sigma _{13}(\xi _{1})}e^{\sigma _{32}(\xi _{2})}+1\otimes
E_{13} \\
-\xi _{2}H_{13}^{\perp }\otimes E_{12}e^{-\sigma _{13}(\xi _{1})}
\end{array}
\right) (1\otimes 1+\xi _{2}C(\xi _{1}))^{-1};
\\[1ex]
\Delta _{\wp }(E_{23}) = & \left(
\begin{array}{c}
E_{23}\otimes e^{-\sigma _{13}(\xi _{1})}-\xi _{2}H_{13}^{\perp }\otimes
H_{23} \\
-\xi _{2}(H_{13}^{\perp })^{2}\otimes e^{-\sigma _{13}(\xi
_{1})}+H_{13}^{\perp }\otimes 1
\end{array}
\right) (1\otimes 1+\xi _{2}C(\xi _{1}))^{-1}- \\[1ex]
&   +\xi _{2}(H_{13}^{\perp }(H_{13}^{\perp }-1)\otimes 1)(1\otimes 1+\xi
_{2}C(\xi _{1}))^{-2}+ \\[1ex]
&   +1\otimes E_{23};
\\[1ex]
\Delta _{\wp }(E_{32}) = & E_{32}\otimes e^{\sigma _{32}(\xi
_{2})}+1\otimes E_{32}+ \\[1ex]
&   +\xi _{1}(\xi _{2}E_{32}+2e^{\sigma _{32}(\xi _{2})}H_{13}^{\perp
})\otimes E_{12}e^{-\sigma _{13}(\xi _{1})}+ \\[1ex]
&   +\xi _{1}^{2}(\xi _{2}^{2}E_{32}+e^{\sigma _{32}(\xi
_{2})}H_{13}^{\perp })H_{13}^{\perp }\otimes (E_{12})^{2}(e^{\sigma
_{13}(\xi _{1})}e^{\sigma _{32}(\xi _{2})}e^{\sigma _{13}(\xi _{1})})^{-1};
\end{array}
\end{equation}
\begin{eqnarray}
\Delta _{\wp }(E_{21})  =  \left\{ E_{21}\otimes e^{-\sigma _{13}(\xi
_{1})}+\xi _{2}H_{13}^{\perp }\otimes E_{31}+\xi _{1}\xi _{2}(H_{13}^{\perp
})^{2}\otimes H_{13}e^{-\sigma _{13}(\xi _{1})}+\right. \nonumber\\[2ex]
+\left(
\begin{array}{c}
\xi _{1}\xi _{2}(H_{13}^{\perp })^{2}- \\
-\xi _{1}\xi _{2}(H_{13}^{\perp })^{3}+\xi _{1}H_{13}^{\perp }E_{23}
\end{array}
\right) \otimes (e^{-\sigma _{13}(\xi _{1})}-e^{-2\sigma _{13}(\xi _{1})})+
\nonumber\\[2ex]
+\xi _{1}H_{23}^{\perp }E_{23}\otimes e^{-\sigma _{13}(\xi _{1})}-\xi
_{1}E_{23}\otimes H_{12}e^{-\sigma _{13}(\xi _{1})}+ \nonumber\\[2ex]
+\xi _{1}^{2}\xi _{2}E_{23}H_{13}^{\perp }\otimes e^{-\sigma _{13}(\xi
_{1})}e^{-\sigma _{32}(\xi _{2})}e^{-\sigma _{13}(\xi _{1})}H_{23}E_{12}+
\nonumber\\[2ex]
\left. +\xi _{1}^{2}\xi _{2}E_{23}(H_{13}^{\perp })^{2}\otimes
E_{12}e^{-\sigma _{13}(\xi _{1})}e^{-\sigma _{32}(\xi _{2})}e^{-2\sigma
_{13}(\xi _{1})}\right\} \times  \nonumber\\[2ex]
\times (1\otimes 1+\xi _{2}C(\xi _{1}))^{-1}+ \nonumber\\
+\left\{ \xi _{1}^{2}\xi _{2}E_{23}(H_{13}^{\perp }-(H_{13}^{\perp
})^{2})\otimes E_{12}e^{-\sigma _{13}(\xi _{1})}e^{-\sigma _{32}(\xi
_{2})}e^{-\sigma _{13}(\xi _{1})}-\right.  \nonumber \\[2ex]
   \left.
-\xi _{1}H_{13}^{\perp }E_{23}\otimes e^{-\sigma _{13}(\xi
_{1})}\right\} (1\otimes 1+\xi _{2}C(\xi _{1}))^{-2}+ \nonumber \\[2ex]
   +\xi _{1}\xi _{2}(1\otimes 1+\xi _{2}C(\xi _{1}))^{-1}\times \nonumber \\[2ex]
   \times \left\{ H_{12}^{\perp }(H_{13}^{\perp }-(H_{13}^{\perp
})^{2})\otimes e^{-\sigma _{32}(\xi _{2})}e^{-\sigma _{13}(\xi
_{1})}\right\} - \nonumber \\
-\xi _{1}H_{23}^{\perp }\otimes E_{23}e^{-\sigma _{13}(\xi _{1})}+ \nonumber \\[2ex]
  +\xi _{1}\xi _{2}H_{12}^{\perp }H_{13}^{\perp }\otimes
(H_{23}-1)e^{-\sigma _{32}(\xi _{2})}e^{-\sigma _{13}(\xi _{1})}+ \nonumber \\
   +(\xi _{1}\xi _{2}H_{12}^{\perp }(H_{13}^{\perp })^{2}-\xi
_{1}H_{12}^{\perp }E_{23})\otimes e^{-\sigma _{13}(\xi _{1})}e^{-\sigma
_{32}(\xi _{2})}e^{-\sigma _{13}(\xi _{1})}- \nonumber \\[2ex]
   -\xi _{1}^{2}\xi _{2}H_{12}^{\perp }H_{13}^{\perp }\otimes
E_{23}E_{12}e^{-\sigma _{13}(\xi _{1})}e^{-\sigma _{32}(\xi _{2})}e^{-\sigma
_{13}(\xi _{1})}- \nonumber \\[2ex]
  -\xi _{1}^{2}E_{23}\otimes e^{-\sigma _{13}(\xi _{1})}e^{-\sigma
_{32}(\xi _{2})}e^{-\sigma _{13}(\xi _{1})}E_{12}E_{23}- \nonumber \\[2ex]
   -\xi _{1}^{2}E_{23}^{2}\otimes E_{12}e^{-2\sigma _{13}(\xi
_{1})}(1\otimes 1+\xi _{2}C(\xi _{1}))^{-1}(1\otimes e^{-\sigma _{32}(\xi
_{2})}e^{-\sigma _{13}(\xi _{1})}) \nonumber \\[2ex]
   +1\otimes E_{21};
\end{eqnarray}
\begin{eqnarray}
\Delta _{\wp }(E_{31})  = \qquad   E_{31}\otimes e^{-\sigma _{13}(\xi
_{1})}-\left\{ \xi _{1}E_{21}+\xi _{1}^{2}(H_{23}^{\perp }-1)E_{23}-\xi
_{1}\xi _{2}H_{13}^{\perp }E_{31}+\right.  \nonumber \\[2ex]
 \left. +\xi _{1}^{2}\xi _{2}(H_{23}^{\perp }-1)H_{12}^{\perp
}H_{13}^{\perp }\right\} \otimes E_{12}e^{-\sigma _{13}(\xi _{1})}e^{-\sigma
_{32}(\xi _{2})}e^{-\sigma _{13}(\xi _{1})}+ \nonumber \\[2ex]
 +\xi _{1}H_{23}^{\perp }\otimes H_{13}e^{-\sigma _{13}(\xi _{1})} \nonumber \\%
[2ex]
 +\xi _{1}H_{12}^{\perp }H_{13}^{\perp }\otimes (e^{-\sigma _{13}(\xi
_{1})}-e^{-\sigma _{32}(\xi _{2})}e^{-\sigma _{13}(\xi _{1})})+ \nonumber \\
 +\xi _{1}(H_{23}^{\perp }-(H_{23}^{\perp })^{2})\otimes (e^{-\sigma
_{13}(\xi _{1})}-e^{-2\sigma _{13}(\xi _{1})})+ \nonumber \\[2ex]
 +(\xi _{1}^{2}E_{23}+\xi _{1}^{2}\xi _{2}H_{12}^{\perp }H_{13}^{\perp
})\otimes H_{13}E_{12}e^{-\sigma _{13}(\xi _{1})}e^{-\sigma _{32}(\xi
_{2})}e^{-\sigma _{13}(\xi _{1})}+ \nonumber \\[2ex]
 +\left\{ 2\xi _{1}^{2}H_{13}^{\perp }E_{23}-\xi _{1}^{2}\xi
_{2}H_{12}^{\perp }H_{13}^{\perp }+\xi _{1}^{2}\xi _{2}(H_{12}^{\perp
})^{2}H_{13}^{\perp }+\right.  \nonumber \\[2ex]
 \left. +2\xi _{1}^{2}\xi _{2}H_{12}^{\perp }(H_{13}^{\perp
})^{2}\right\} \otimes E_{12}e^{-2\sigma _{13}(\xi _{1})}e^{-\sigma
_{32}(\xi _{2})}e^{-\sigma _{13}(\xi _{1})}+ \nonumber \\[2ex]
 +\left\{ 2\xi _{1}^{2}(H_{12}^{\perp }-1)E_{23}-\xi _{1}^{2}\xi
_{2}H_{12}^{\perp }H_{13}^{\perp }+\right.  \nonumber \\[2ex]
 \left. +\xi _{1}^{2}\xi _{2}(H_{12}^{\perp })^{2}H_{13}^{\perp
}\right\} \otimes E_{12}e^{-\sigma _{13}(\xi _{1})}e^{-\sigma _{32}(\xi
_{2})}e^{-2\sigma _{13}(\xi _{1})}+ \nonumber \\
 +\left\{ 2\xi _{1}^{3}\xi _{2}(H_{12}^{\perp }-1)(H_{13}^{\perp
}+1)E_{23}+\xi _{1}^{3}\xi _{2}^{2}(H_{12}^{\perp })^{2}(H_{13}^{\perp
})^{2}-\right.  \nonumber \\[2ex]
 \left. -\xi _{1}^{3}\xi _{2}^{2}H_{12}^{\perp }(H_{13}^{\perp
})^{2}\right\} \otimes (E_{12}e^{-\sigma _{13}(\xi _{1})}e^{-\sigma
_{32}(\xi _{2})}e^{-\sigma _{13}(\xi _{1})})^{2}+ \nonumber \\[2ex]
 +\xi _{1}(E_{23}\otimes e^{-\sigma _{13}(\xi _{1})}E_{32})(1\otimes
1+\xi _{2}C(\xi _{1}))^{-1}+ \nonumber \\[2ex]
 +\xi _{1}^{2}((H_{13}^{\perp }+1)E_{23}\otimes E_{12}e^{-2\sigma
_{13}(\xi _{1})})(1\otimes 1+\xi _{2}C(\xi _{1}))^{-1}- \nonumber \\[2ex]
 -\xi _{1}^{2}(1\otimes 1+\xi _{2}C(\xi _{1}))^{-1}\times  \nonumber \\[2ex]
 \times (H_{13}^{\perp }E_{23}\otimes E_{12}e^{-\sigma _{13}(\xi
_{1})}e^{-\sigma _{32}(\xi _{2})}e^{-\sigma _{13}(\xi _{1})})+ \nonumber \\
 +\xi _{1}^{3}(E_{23}^{2}\otimes E_{12}^{2}e^{-2\sigma _{13}(\xi
_{1})})(1\otimes 1+\xi _{2}C(\xi _{1}))^{-1}\times  \nonumber \\[2ex]
 \times (1\otimes e^{-\sigma _{32}(\xi _{2})}e^{-\sigma _{13}(\xi
_{1})})(1\otimes 1+\xi _{2}C(\xi _{1}))\times  \nonumber \\
 \times (1\otimes e^{-\sigma _{32}(\xi _{2})}e^{-\sigma _{13}(\xi
_{1})})+ \nonumber \\
 +1\otimes E_{31};
\label{copr-par}
\end{eqnarray}
where
\[
C(\xi_{1})=1\otimes E_{32}+\xi_{1}H_{13}^{\perp }\otimes
E_{12}e^{-\sigma _{13}(\xi_{1})}.
\]
The parabolic twist $\mathcal{F}_{\wp }$ can be supplied with two
natural parameters corresponding to two Jordanian-like factors:
\begin{equation}
\begin{array}{l}
\mathcal{F}_{\wp }(\xi ,\zeta )= \\
\exp (H_{13}^{\perp }\otimes (2\sigma _{13}(\xi )+\sigma _{32}(\zeta )))\exp
(-\xi E_{23}\otimes E_{12}e^{\sigma _{13}(\xi )})\exp (H_{23}\otimes \sigma
_{13}(\xi ))
\end{array}
\label{par-param}
\end{equation}

If in the universal $\mathcal{R%
}$-matrix for $U_{\wp }(\frak{P};\xi ,\zeta )$,
\begin{equation}
 \begin{array}{c}
\mathcal{R}_{\wp }(\xi ,\zeta )=\left( \mathcal{F}_{\wp }(\xi ,\zeta
)\right) _{21}\left( \mathcal{F}_{\wp }(\xi ,\zeta )\right) ^{-1}= \\
\exp ((2\sigma _{13}(\xi )+\sigma _{32}(\zeta ))\otimes H_{13}^{%
\scriptscriptstyle\perp })\exp (-\xi E_{12}e^{\sigma _{13}(\xi )}\otimes
E_{23})\exp (\sigma _{13}(\xi )\otimes H_{23})\times  \\
\times \exp (-H_{23}\otimes \sigma _{13}(\xi ))\exp (\xi E_{23}\otimes
E_{12}e^{\sigma _{13}(\xi )})\exp (-H_{13}^{\scriptscriptstyle\perp }\otimes
(2\sigma _{13}(\xi )+\sigma _{32}(\zeta )).
\end{array}
\label{rrmat-param}
\end{equation}
the parameters are chosen to be proportional ($\zeta =\eta \xi $) the
expression (\ref{rrmat-param}) can be considered as a quantization of the
classical $r$-matrix
\begin{equation*}
r_{\wp }(\eta )=H_{23}^{\perp }\wedge E_{13}+E_{12}\wedge E_{23}+\eta
H_{13}^{\perp }\wedge E_{32}.
\end{equation*}

\section{Quantum twists for quasi-Frobenius subalgebras in $\frak{sl_3}$.}

In what follows we define quantum deformations ($q$-twists) for the twists
constructed in the previous section so that
the diagram (\ref{q-tw}) commutes. For the quantum algebra $U_{q}(\frak{sl}_{3})$
the generators will be denoted by the small letters and we shall use the following
defining relations:
\begin{equation*}
\begin{array}{lcl}
\lbrack h_{ij},e_{kl}]=(\delta _{ik}+\delta _{jl}-\delta _{il}-\delta
_{jk})\,\,\,e_{kl}, &  & [e_{12},e_{32}]=[e_{21},e_{23}]=0
\end{array}
\end{equation*}
\begin{equation*}
\begin{array}{lclclcl}
\lbrack e_{12},e_{21}] & = & \frac{q^{h_{12}}-q^{-h_{12}}}{q-q^{-1}}, &  &
[e_{23},e_{32}] & = & \frac{q^{h_{23}}-q^{-h_{23}}}{q-q^{-1}} \\[2ex]
e_{13}e_{12} & = & q^{-1}~e_{12}e_{13}, &  & e_{13}e_{23} & = &
q~e_{23}e_{13} \\[2ex]
e_{21}e_{31} & = & q\;e_{31}e_{21}, &  & e_{32}e_{31} & = &
q^{-1}\;e_{31}e_{32}
\end{array}
\end{equation*}
where the composite root generators $e_{13}$ and $e_{31}$ are defined as
follows
\begin{equation*}
\begin{array}{lcl}
e_{13}:=e_{12}e_{23}-q^{-1}~e_{23}e_{12}, &  &
e_{31}:=e_{32}e_{21}-q~e_{21}e_{32}
\end{array}
\end{equation*}
and the coproduct is fixed by its values on the Chevalley generators
\begin{equation*}
\begin{array}{lcl}
\Delta (h_{ij})=h_{ij}\otimes 1+1\otimes h_{ij}, &  & \Delta
(e_{i,i+1})=q^{-h_{i,i+1}}\otimes e_{i,i+1}+e_{i,i+1}\otimes 1 \\[2ex]
&  & \Delta (e_{i+1,i})=e_{i+1,i}\otimes q^{h_{i,i+1}}+1\otimes e_{i+1,i}.
\end{array}
\end{equation*}

$q-$twists are defined for the deformed carrier Hopf
subalgebras in $U_{q}(\mathfrak{sl}_{3})$. We consider these
carrier subalgebras as $q-$quantization of the classical
quasi-Frobenius subalgebras in $\mathfrak{sl}_{3}$.

\subsection{Abelian two dimensional subalgebras}

\subsubsection{$q$-deformation for $\frak{h}$}

As far as in this case the carrier in $U_{q}(\frak{sl}_{3})$
\begin{equation*}
\begin{array}{lcl}
H & = & \left( {\
\begin{array}{ccc}
\ast  & 0 & 0 \\
0 & \ast  & 0 \\
0 & 0 & \ast
\end{array}
}\right)
\end{array}
\end{equation*}
is undeformed,
\begin{equation*}
\Delta (h_{ij})=h_{ij}\otimes 1+1\otimes h_{ij},
\end{equation*}
the corresponding Abelian twist
\begin{equation}
\mathcal{F}_{\mathcal{R}}=\exp (t\xi ^{ij,kl}h_{ij}\otimes h_{kl})\quad \quad \xi
^{ij,kl}\in \mathbb{C}.  \label{ab-h-3}
\end{equation}
can be taken $q-$independent (with the limit  $\exp (t\xi ^{ij,kl}H_{ij}\otimes H_{kl})$).

\subsubsection{$q-$quantization of $\frak{h}^{(1)}$}

By definition
\begin{equation*}
\begin{array}{lcl}
\frak{h}^{(1)} & = & \ast \left( {\
\begin{array}{ccc}
1 & 0 & 0 \\
0 & 1 & 0 \\
0 & 0 & -2
\end{array}
}\right) +\left( {\
\begin{array}{ccc}
0 & \ast  & 0 \\
0 & 0 & 0 \\
0 & 0 & 0
\end{array}
}\right)
\end{array}
\end{equation*}
and $U_{q}(\frak{h}^{(1)})$ is the following Hopf subalgebra of $U_{q}(\frak{%
sl}_{3})$
\begin{equation*}
\begin{array}{lcl}
\Delta (e_{12})=q^{-h_{12}}\otimes e_{12}+e_{12}\otimes 1, &  & \Delta
(h_{12}^{\perp })=h_{12}^{\perp }\otimes 1+1\otimes h_{12}^{\perp }
\end{array}
\end{equation*}
where $h_{12}^{\perp }=\frac{1}{3}(e_{11}+e_{22})-\frac{2}{3}e_{33}$. The
next step is to define a contraction of the algebra $U_{q}(\frak{sl}_{3})$
leading to the deformation of $U(\frak{h}^{(1)})$ defined by the
twisting element (see (\ref{ab-h-3}))
\begin{equation*}
\mathcal{F}=\exp (\xi ~H_{\lambda }^{\perp }\otimes E_{\lambda }).
\end{equation*}
To find such limiting procedure we introduce a family of Hopf algebras
equivalent to $U_{q}(\frak{h}^{(1)})$ . This family is obtained by applying
to $U_{q}(\frak{h}^{(1)})$ the similarity transformation defined by the
coboundary twist. To fix its form we use the notations
$$
\exp _{q}\left( x \right)=\sum_{n \geq 0}\frac{x^n}{\left( n \right)_{q} !}
$$
where
\begin{equation*}
\begin{array}{lcl}
(n)_{q}=\frac{q^{n}-1}{q-1}, &  & (n)_{q}!=(1)_{q}(2)_{q}\cdots (n)_{q},
\end{array}
\end{equation*}
Now put $W=\exp _{q^{2}}(ts^{-1}\;e_{12})$. The necessary coboundary twist
$\mathcal{J}\left( s,t\right)$ is
\begin{equation}
{\
\begin{array}{r}
\mathcal{J}\left( s,t\right) :=(W\otimes W)\Delta (W^{-1}).
\end{array}
}  \label{asympt}
\end{equation}

According to the Heine formula \cite{Kac}
\begin{equation*}
\begin{array}{r}
(1-t\;x)_{q}^{(-\alpha )}:=1+\sum_{n\geq 1}t^{n}\frac{(\alpha )_{q}\cdots
(\alpha +n-1)_{q}}{(n)_{q}!}x^{n}\\
=\exp _{q}(\frac{t}{1-q}~x)\exp _{q^{-1}}(-%
\frac{q^{\alpha }t}{1-q}~x)
\end{array}
\end{equation*}
Thus the twisting element $\mathcal{J}\left( s,t\right) $ can be simplified,
\begin{equation*}
\mathcal{J}\left( s,t\right) =1\otimes 1+\sum_{n\geq
1}(s^{-1}(1-q^{2}))^{n}t^{n}~\frac{(-\frac{1}{2}h_{12})_{q^{2}}\cdots (-%
\frac{1}{2}h_{12}+n-1)_{q^{2}}\otimes e_{12}^{n}}{(n)_{q^{2}}}.
\end{equation*}
If we require that
\begin{equation}
q\equiv 1+s^{2}t~{\rm mod}\left(s^{3}t \right)  \label{q}
\end{equation}
then $\mathcal{J}\left( s,t\right) $ contains only positive degrees of $s$
and moreover
\begin{equation}
\mathcal{J}\left( s,t\right) =1\otimes 1~{\rm mod}\left(st \right).
\label{cob-approx}
\end{equation}
Now consider the element
\begin{equation}
F_{q}:=(W\otimes W)q^{s^{-1}\cdot h_{12}^{\perp }\otimes h_{12}}\Delta
(W^{-1})=\left( \mathrm{Ad}(W\otimes W)\circ (q^{s^{-1}\cdot h_{12}^{\perp
}\otimes h_{12}})\right) \mathcal{J}\left( s,t\right) .  \label{q-twist-1}
\end{equation}
It defines a twist for $U_{q}(\frak{sl}_{3})$ because it satisfies the
Drinfeld equation being equivalent to the Abelian twist $q^{s^{-1}\cdot
h_{12}^{\perp }\otimes h_{12}}$.

Let us check that
\begin{equation}
F_{q}\equiv \exp (-2t^{2}~h_{12}^{\perp }
\otimes e_{12}){\rm mod}\left(st \right).
\label{q-tw-1-approx}
\end{equation}
Notice that the multiplier $t^{2}$ is required as far as $q$-dependent
terms must not contribute to the twist in the limit $s \rightarrow 0$.
Using the Heine's formula we can calculate explicitly
\begin{equation*}
\mathrm{Ad}(W\otimes W)\circ (q^{s^{-1}\cdot ~h_{12}^{\perp }\otimes
h_{12}})=(1\otimes \mathrm{Ad}(W)\circ (q^{h_{12}}))^{s^{-1}\cdot
~h_{12}^{\perp }\otimes 1}
\end{equation*}
and
\begin{equation*}
\mathrm{Ad}(W)\circ (q^{h_{12}})=\frac{1}{1-(1-q^{2})s^{-1}t\cdot e_{12}}%
q^{h_{12}}\equiv (1-2st^{2}~e_{12}){\rm mod}\left(s^{2}t \right).
\end{equation*}
This together with (\ref{cob-approx}) and (\ref{q-twist-1}) proves (\ref
{q-tw-1-approx}). In the limit $s\longrightarrow 0$ we get the special case
of the general twisting element $\mathcal{F}_{\mathcal{R}}=\exp (\xi
~H_{\lambda }^{\perp }\otimes E_{\lambda })$ (see (\ref{ab-h1})).

\subsubsection{$q-$quantization of $\frak{h}^{(0,1)}$}

By definition
\begin{equation*}
\begin{array}{lcl}
\frak{h}^{(0,1)} & = & \left( {\
\begin{array}{ccc}
0 & \ast  & \ast  \\
0 & 0 & 0 \\
0 & 0 & 0
\end{array}
}\right)
\end{array}
\end{equation*}
To find the necessary quantized carrier $U_{q}(\frak{h}^{(0,1)})$ let us
simplify the form of the corresponding Hopf subalgebra in $U_{q}(\frak{sl}%
_{3})$ . There in particular we have
\begin{equation*}
\Delta (e_{13})=q^{-h_{13}}\otimes e_{13}+e_{13}\otimes
1+(1-q^{-2})~q^{-h_{12}}e_{23}\otimes e_{12}.
\end{equation*}
Let us perform the twist transformation $U_{q}(\frak{sl}_{3})\longrightarrow
U_{q,R_{1}}(\frak{sl}_{3})$ by applying the  $R-$matrix factor twisting
element
\begin{equation*}
R_{1}=\exp _{q^{2}}(-(q-q^{-1})~e_{21}\otimes e_{12}).
\end{equation*}
In the twisted algebra $U_{q,R_{1}}(\frak{sl}_{3})$ the coalgebra of $%
U_{q}(\frak{h}^{(0,1)})$ is generated by
\begin{equation*}
\begin{array}{lcl}
\Delta (e_{12})=q^{h_{12}}\otimes e_{12}+e_{12}\otimes 1, &  & \Delta
(e_{13})=q^{-h_{13}}\otimes e_{13}+e_{13}\otimes 1.
\end{array}
\end{equation*}

To define the quantum analogue of the twist $\mathcal{F}_{\mathcal{R}}=~\exp
(\xi ~E_{12}\otimes ~E_{13})$ -- the special case of the general expression (%
\ref{ab-h01-11}) with $\mu =0$ -- we consider the following $q$-twist:
\begin{equation*}
F_{q}=(W\otimes W)R_{1}\Delta (W^{-1}),\qquad W=\exp
_{q^{2}}(s^{-1}t~q^{-h_{12}}e_{12})\exp _{q^{2}}(s^{-1}t~e_{13})
\end{equation*}
(Notice that $[q^{-h_{12}}e_{12},e_{13}]=0$).\newline
Explicitly,
\begin{equation*}
\begin{array}{r}
F_{q}=\left( \exp _{q^{2}}(s^{-1}t~q^{-h_{12}}e_{12})\otimes \exp
_{q^{2}}(s^{-1}t~e_{13})\right) \exp _{q^{-2}}(-ts^{-1}~q^{-h_{13}}\otimes
e_{13})\cdot  \\[2ex]
\cdot \exp _{q^{-2}}(-s^{-1}t~q^{-h_{12}}e_{12}\otimes
q^{-h_{12}})R_{1}.
\end{array}
\end{equation*}
Using the relation
\begin{equation*}
\mathrm{Ad}(\exp _{q^{2}}(s^{-1}t~q^{-h_{12}}e_{12}))\circ
(q^{-h_{13}})=(1-(1-q^{2})s^{-1}t\cdot q^{-h_{12}}e_{12})_{q^{2}}^{(\frac{1}{%
2})}q^{-h_{13}}
\end{equation*}
and assuming that $q\equiv 1+s^{2}t~{\rm mod}\left(s^{3}t \right)$ we can check that
\begin{equation*}
F_{q}\equiv \exp (t^{3}~e_{12}\otimes e_{13})~{\rm mod}\left(st \right).
\end{equation*}
Another possible $q$-twist corresponding to the same
bialgebraic structure looks like
\begin{equation*}
F_{q}^{\prime}=\exp_{q^{2}}(t\mathop{}q^{-h_{12}}e_{12}\otimes
e_{13})q^{-h_{13}^{\perp}\otimes h_{12}^{\perp}}R_{1}.
\end{equation*}

\subsubsection{$q-$quantization of $\frak{h}^{(1,1)}$  by embedding into $%
U_{q}(A_{2}^{(2)})$}

By definition
\begin{equation*}
\begin{array}{lcc}
\frak{h}^{(1,1)} & = & \ast \left( {\
\begin{array}{ccc}
0 & 1 & 0 \\
0 & 0 & 1 \\
0 & 0 & 0
\end{array}
}\right) +\left( {\
\begin{array}{ccc}
0 & 0 & \ast  \\
0 & 0 & 0 \\
0 & 0 & 0
\end{array}
}\right)
\end{array}
\end{equation*}
The $q$-quantization of the classical $r$-matrix $\rho=E_{13}\wedge
(E_{12}+E_{23})$ can be related to the twist for
$U_{q}(\hat{\mathfrak{sl}}_{3})$. Following \cite{Sam} we define here
the $q$-twist in the root generators notation:
$$
\begin{array}{l}
{\cal F}_{q}=\\[2ex]
\exp_{q^{2}}(\frac 12 t\mathop{}\hat{e}_{\delta-\alpha-\beta}
\otimes\hat{e}_{-\beta})\exp_{q^{2}}(t\mathop{}
\hat{e}_{\delta-\alpha-\beta}\otimes\hat{e}_{-\alpha})
\exp_{q^{2}}(-\frac 12 qt\mathop{}\hat{e}_{\delta-\beta}\otimes
\hat{e}_{-\alpha-\beta})\cdot\\[2ex]
\hphantom{aaaaaaaaaaaaaaaaaaaaaaaaaaaaaaaaaaaaa}\exp_{q^{2}}(-\frac
12 (q-q^{-1})\mathop{}\hat{e}_{\alpha} \otimes
\hat{e}_{-\beta}){\cal K}
\end{array}
$$
where
$$
 {\cal K}=q^{\frac 49 h_{\alpha}\otimes h_{\alpha}+ \frac 29
h_{\alpha}\otimes h_{\beta}+ \frac 59 h_{\beta}\otimes
h_{\alpha}+\frac 79 h_{\beta}\otimes h_{\beta}},
$$
$$
\begin{array}{lll}
\hat{e}_{-\alpha}=q^{\frac
12\mathop{}h_{\beta}^{\perp}}e_{-\alpha},&\hat{e}_{-\beta}=q^{-\frac
12\mathop{}h_{\beta}^{\perp}}e_{-\beta},&
\hat{e}_{\delta-\alpha-\beta}=q^{-\frac
12\mathop{}h_{\alpha+\beta}^{\perp}}
e_{\delta-\alpha-\beta},\\[2ex]
\hat{e}_{\alpha}=e_{\alpha}q^{\frac
12\mathop{}h_{\alpha+\beta}^{\perp}},&
\hat{e}_{\delta-\beta}^{\prime}=[\hat{e}_{\alpha},\hat{e}_{\delta-\alpha-\beta}],
\end{array}
$$
and
$$
\begin{array}{llll}
h_{\alpha}^{\perp}=\frac 23 h_{\alpha}+\frac 43 h_{\beta},&
h_{\beta}^{\perp}=\frac 43 h_{\alpha}+\frac 23 h_{\beta},&
h_{\alpha+\beta}^{\perp}=h_{\beta}^{\perp}-h_{\alpha}^{\perp},&
\hat{e}_{-\beta}=q^{-\frac
12\mathop{}h_{\beta}^{\perp}}e_{-\beta}.
\end{array}
$$

In the limit $s\rightarrow 0$ and assuming that $q=1+st$, we come
to the following twist for
$U(\widehat{\mathfrak{sl}}_{3})$:
$$
\begin{array}{l}
{\cal F}_{1}=\\[2ex]
\exp(\frac 12 t\mathop{}E_{\delta-\alpha-\beta} \otimes
E_{-\beta})\exp(t\mathop{} E_{\delta-\alpha-\beta}\otimes
E_{-\alpha}) \exp (-\frac 12 t\mathop{}E_{\delta-\beta}\otimes
E_{-\alpha-\beta})\\[2ex]
\end{array}
$$
which can be considered as an affinization of a twist quantizing
the $r$-matrix
$$
\widehat{\rho}=-\frac 12 E_{\delta-\alpha-\beta}\wedge
E_{-\beta}+\frac 12 E_{\delta-\beta}\wedge
E_{-\alpha-\beta}-E_{\delta-\alpha-\beta}\wedge E_{-\alpha}.
$$

On the other hand we
can consider the quantum twisted affine Hopf algebra
$U_{q}(A_{2}^{(2)})$ (see also \cite {TL}) i.e. the Drinfeld-Jimbo
quantization of the Cartan matrix:
\begin{equation*}
A=\left( {\
\begin{array}{rr}
2 & -1 \\[2ex]
-4 & 2
\end{array}
}\right) =DB=\left( {\
\begin{array}{rr}
\frac{1}{2} & 0 \\[2ex]
0 & 2
\end{array}
}\right) \left( {\
\begin{array}{rr}
4 & -2 \\[2ex]
-2 & 1
\end{array}
}\right) .
\end{equation*}
In the evaluation representation we have the Kac generators
\begin{equation*}
\begin{array}{rcl}
H_{0}=-2H_{13}, &  & H_{1}=H_{13}, \\[2ex]
E_{0}=\sqrt{2}E_{31}u, &  & E_{1}=(E_{12}+E_{23}), \\[2ex]
F_{0}=\sqrt{2}E_{13}u^{-1}, &  & F_{1}=(E_{21}+E_{32}).
\end{array}
\end{equation*}
The Drinfeld-Jimbo
quantization $U_{q}(A_{2}^{(2)})$ is defined by the relations
\begin{equation*}
\begin{array}{lcl}
\lbrack h_{\alpha },e_{\delta -2\alpha }]=-4e_{\delta -2\alpha }, &  &
[h_{\alpha },e_{\alpha }]=2e_{\alpha }, \\[2ex]
\lbrack e_{\delta -2\alpha },e_{-\delta +2\alpha }]=\frac{q^{-h_{\alpha
}}-q^{h_{\alpha }}}{q-q^{-1}}, &  & [e_{\alpha },e_{-\alpha }]=\frac{q^{%
\frac{1}{2}h_{\alpha }}-q^{-\frac{1}{2}h_{\alpha }}}{q-q^{-1}},
\end{array}
\end{equation*}
\begin{equation*}
\begin{array}{l}
\Delta (e_{\delta -2\alpha })=q^{h_{\alpha }}\otimes e_{\delta -2\alpha
}+e_{\delta -2\alpha }\otimes 1, \qquad   \Delta (e_{\alpha })=q^{-\frac{1}{2}%
h_{\alpha }}\otimes e_{\alpha }+e_{\alpha }\otimes 1, \\[2ex]
\Delta (e_{-\delta +2\alpha })=e_{-\delta +2\alpha }\otimes q^{-h_{\alpha
}}+1\otimes e_{-\delta +2\alpha }, \,  \Delta (e_{-\alpha })=e_{-\alpha
}\otimes q^{\frac{1}{2}h_{\alpha }}+1\otimes e_{-\alpha },
\end{array}
\end{equation*}
plus the Serre relations of the form:
\begin{equation*}
\begin{array}{lcl}
\lbrack \lbrack e_{\alpha },e_{\delta -2\alpha }]_{q},e_{\delta -2\alpha
}]_{q}=0, &  & [e_{\alpha },[e_{\alpha },[e_{\alpha },[e_{\alpha
},[e_{\alpha },e_{\delta -2\alpha }]_{q}]_{q}]_{q}]_{q}]_{q}=0
\end{array}
\end{equation*}
where
\begin{equation*}
\lbrack e_{i},e_{j}]_{q}:=xy-q^{b_{ij}}yx.
\end{equation*}
Let us fix the normal ordering,
\begin{equation*}
\alpha \prec \delta +2\alpha \prec \delta +\alpha \prec 3\delta +2\alpha
\prec 2\delta +\alpha \prec \delta \prec 2\delta -\alpha \prec 3\delta
-2\alpha \prec \delta -\alpha \prec \delta -2\alpha,
\end{equation*}
and define the corresponding ordering on the set of Chevalley generators
in $U(A_{2}^{(2)})$:
\begin{equation*}
E_{1}\prec E_{8}\prec E_{5}\prec E_{9}\prec E_{7}\prec E_{3}\prec E_{4}\prec
E_{6}\prec E_{2}\prec E_{0}
\end{equation*}
where
\begin{equation*}
\begin{array}{lcl}
E_{2}=\sqrt{2}(E_{21}-E_{32})u, &  & E_{3}=\sqrt{2}(H_{12}-H_{23})u, \\[2ex]
E_{4}=-3\sqrt{2}(E_{12}-E_{23})u, &  & E_{5}=6\sqrt{2}E_{13}u.
\end{array}
\end{equation*}
Define the twisting element
\begin{equation*}
F_{q}=(W\otimes W)\Delta (W^{-1}),\qquad W=\exp _{q}(ts^{-1}~e_{\alpha
})\exp _{q^{4}}(ts^{-1}~e_{-\delta +2\alpha }).
\end{equation*}
Explicitly,
\begin{eqnarray*}
F_{q} &=&\exp _{q^{4}}(ts^{-1}~1\otimes e_{-\delta +2\alpha })\exp
_{q^{-4}}(-ts^{-1}~K\otimes e_{-\delta +2\alpha })\cdot  \\
&&\cdot \exp _{q}(ts^{-1}~1\otimes e_{\alpha })\exp _{q^{-1}}(-ts^{-1}~q^{-%
\frac{1}{2}h}\otimes e_{\alpha }),
\end{eqnarray*}
where
\begin{equation*}
K:=(1-(1-q)s^{-1}t~e_{\alpha })_{q}^{(2)}q^{-h_{\alpha }}.
\end{equation*}
Imposing the relation $q\equiv 1+s^{2}t~{\rm mod}\left(s^{3}t \right)$, we can check
that
\begin{equation*}
F_{q}\equiv \exp (2t^{3}~e_{\alpha }\otimes e_{-\delta +2\alpha })~{\rm mod}%
\left(st \right).
\end{equation*}
In the limit $s\longrightarrow 0$ we come to the twisting element
\begin{equation*}
\exp (2\sqrt{2}t^{3}u~(E_{12}+E_{23})\otimes E_{13})
\end{equation*}
that is the other special case of the general solution (\ref{ab-h01-11}),
this time with $\mu =1$:
\begin{equation*}
\mathcal{F}_{\mathcal{R}}=\exp (\xi (E_{12}+E_{23})\otimes E_{13}).
\end{equation*}

\subsection{Non-Abelian two-dimensional subalgebras}

We have three types of nonequivalent non-Abelian quasi-Frobenius Lie
subalgebras in $\frak{sl}_{3}$
\begin{equation*}
\begin{array}{lcl}
\frak{b}^{(0)}, & \frak{b}_{\lambda }, & \frak{b}^{(1)}.
\end{array}
\end{equation*}

\subsubsection{$q-$quantization of $\frak{b}^{(0)}$}

In the case $\frak{b}^{(0)}$
\begin{equation}
\begin{array}{lcc}
\frak{b}^{(0)} & = & \ast \left( {\
\begin{array}{ccc}
1 & 0 & 0 \\
0 & 0 & 0 \\
0 & 0 & -1
\end{array}
}\right) +\ast \left( {\
\begin{array}{ccc}
0 & 1 & 0 \\
0 & 0 & 1 \\
0 & 0 & 0
\end{array}
}\right)
\end{array}
\end{equation}
we introduce the quantum twist
\begin{equation*}
F_{q}=(W\otimes W)\Delta (W^{-1}),\qquad W=\exp _{q}(ts^{-1}~e_{\alpha }).
\end{equation*}
Explicitly,
\begin{equation*}
F_{q}=(1-(1-q)s^{-1}t~1\otimes e_{\alpha })_{q}^{-\frac{1}{2}~(h_{\alpha
}\otimes 1)},
\end{equation*}
and put $q\equiv 1+st~{\rm mod}\left(s^{2}t \right)$ . Then in the limit $s\rightarrow 0$
we come to the twist
\begin{equation*}
\mathcal{F}=\exp (H_{13}\otimes \ln (1+t~(E_{12}+E_{23})))
\end{equation*}
(see (\ref{reg-2})).

\subsubsection{$q-$quantization of $\frak{b}_{\protect\lambda }$}

By definition
\begin{equation*}
\frak{b}_{\lambda }=\ast \left( {\
\begin{array}{ccc}
\lambda  & 0 & 0 \\
0 & \lambda -1 & 0 \\
0 & 0 & 1-2\lambda
\end{array}
}\right) +\left( {\
\begin{array}{ccc}
0 & \ast  & 0 \\
0 & 0 & 0 \\
0 & 0 & 0
\end{array}
}\right)
\end{equation*}
The Hopf algebra $U_{q}(\frak{b}_{\lambda })$ is the Hopf subalgebra
containing $e_{12}$ and $h_{23}-\lambda h_{12}^{\perp }$ with the coproducts
\begin{eqnarray*}
\Delta \left( h_{23}-\lambda h_{12}^{\perp }\right)  &=&\left(
h_{23}-\lambda h_{12}^{\perp }\right) \otimes 1+1\otimes \left(
h_{23}-\lambda h_{12}^{\perp }\right) , \\
\Delta (e_{12}) &=&q^{2h_{23}-2\lambda h_{12}^{\perp }}\otimes
e_{12}+e_{12}\otimes 1.
\end{eqnarray*}
Note that we have the embedding $U_{q}(\frak{b}_{\lambda })\hookrightarrow
U_{q,\mathcal{K}_{\lambda }}(\frak{sl}_{3})$ into the twisted algebra $%
U_{q,\mathcal{K}_{\lambda }}(\frak{sl}_{3})$ where the corresponding
twisting element is
\begin{equation*}
\mathcal{K}_{\lambda }=q^{(2\lambda -3)~h_{12}^{\perp }\otimes h_{23}}.
\end{equation*}
Define the quantum twist
\begin{eqnarray*}
F_{q} &=&\left( \exp _{q^{2}}(s^{-1}t~e_{12})\otimes \exp
_{q^{2}}(s^{-1}t~e_{12})\right) \Delta (\exp _{q^{-2}}(-s^{-1}t~e_{12}))= \\
&=&(1-(1-q^{2})s^{-1}t~1\otimes e_{12})_{q^{2}}^{(-(h_{23}-\lambda
~h_{12}^{\perp })\otimes 1)}
\end{eqnarray*}
and put $q\equiv 1+st~{\rm mod}\left(s^{2}t \right)$ . In face of the evaluation
\begin{equation*}
F_{q}\equiv \exp ((-h_{23}+\lambda h_{12}^{\perp })\otimes \ln (1+2t^{2}~e_{12}))~%
{\rm mod}\left(st \right).
\end{equation*}
we see that the desired quantization of the Jordanian twist (\ref{reg-1}) is
obtained.

\subsubsection{$q-$quantization of $\frak{b}^{(1)}$}

This case is given by the $\frak{sl}_{3}$ subalgebra
\begin{equation*}
\frak{b}^{(1)}=\ast \left(
\begin{array}{ccc}
2 & 0 & 0 \\
0 & -1 & 1 \\
0 & 0 & -1
\end{array}
\right) +\left(
\begin{array}{ccc}
0 & 0 & \ast  \\
0 & 0 & 0 \\
0 & 0 & 0
\end{array}
\right) .
\end{equation*}
In $U_{q}(\frak{sl}_{3})$ we have the coproducts:
\begin{equation*}
\begin{array}{l}
\Delta (e_{23})=q^{-h_{23}}\otimes e_{23}+e_{23}\otimes 1, \\
\Delta (e_{13}^{\prime })=q^{-h_{13}}\otimes e_{13}^{\prime }+e_{13}^{\prime
}\otimes 1+(q^{-1}-q)~q^{-h_{23}}e_{12}\otimes e_{23}, \\
e_{13}^{\prime }=e_{12}e_{23}-q~e_{23}e_{12}.
\end{array}
\end{equation*}
Let us twist these coproducts by the $R-$matrix factor
\begin{equation*}
R_{2}=\exp _{q^{2}}(-(q-q^{-1})~e_{32}\otimes e_{23}).
\end{equation*}
This leads to the simplified coalgebra,
\begin{equation*}
\begin{array}{lcl}
\Delta _{R_{2}}(e_{23})=q^{h_{23}}\otimes e_{23}+e_{23}\otimes 1, &  &
\Delta _{R_{2}}(e_{13}^{\prime })=q^{-h_{13}}\otimes e_{13}^{\prime
}+e_{13}^{\prime }\otimes 1.
\end{array}
\end{equation*}
Consider the twisting factor
\begin{equation*}
\mathcal{K}=q^{-h_{13}^{\perp }\otimes h_{23}^{\perp }},\qquad h_{13}^{\perp }=\frac{1%
}{3}(e_{11}+e_{33})-\frac{2}{3}e_{22},
\end{equation*}
then we obtain the coproducts
\begin{eqnarray*}
\Delta _{\mathcal{K}R_{2}}(q^{-h_{23}^{\perp }}e_{23}) &=&q^{-2h_{13}^{\perp }}\otimes
q^{-h_{23}^{\perp }}e_{23}+q^{-h_{23}^{\perp }}e_{23}\otimes 1, \\
\Delta _{\mathcal{K}R_{2}}(e_{13}^{\prime }) &=&q^{-2h_{23}^{\perp }}\otimes
e_{13}^{\prime }+e_{13}^{\prime }\otimes 1.
\end{eqnarray*}
Note that we have the relation
\begin{equation*}
\lbrack q^{-h_{23}^{\perp }}e_{23},e_{13}^{\prime }]=0.
\end{equation*}
$U_{q}(\frak{b}^{(1)})$ is defined as the minimal Hopf subalgebra in $%
U_{q,\mathcal{K}R_{2}}(\frak{sl}_{3})$ containing $e_{13}^{\prime
},e_{23},h_{23}^{\perp }$.

The quantum twist with the necessary limit properties will be constructed in
terms of thus defined Hopf algebra $U_{q}(\frak{b}^{(1)})$ . It contain two
factors. The first one is a coboundary twist of the form
\begin{eqnarray*}
F_{q}^{1} &=&(W\otimes W)\Delta (W^{-1}), \\
W &=&\exp _{q^{-2}}(-s^{-1}t~q^{-h_{23}^{\perp }}e_{23})\exp
_{q^{2}}(s^{-2}t~e_{13}^{\prime }).
\end{eqnarray*}
Explicitly,
\begin{equation*}
\begin{array}{r}
F_{q}^{1}=\exp _{q^{-2}}(-s^{-1}t~1\otimes q^{-h_{23}^{\perp }}e_{23})\exp
_{q^{2}}(s^{-1}t~q^{-2h_{13}^{\perp }}\otimes q^{-h_{23}^{\perp
}}e_{23})\cdot  \\[2ex]
\cdot \exp _{q^{2}}(s^{-2}t~1\otimes e_{13}^{\prime })\exp
_{q^{-2}}(-s^{2}t~q^{-2h_{23}^{\perp }}\otimes e_{13}^{\prime }).
\end{array}
\end{equation*}
If
\begin{equation*}
q\equiv 1+s^{2}t~{\rm mod}\left(s^{3}t \right)
\end{equation*}
then
\begin{equation*}
F_{q}^{1}\equiv \exp (h_{23}^{\perp }\otimes \ln (1+2t^{2}~e_{13}^{\prime
}))~{\rm mod}\left(st \right).
\end{equation*}

In the deformed Hopf algebra $U_{q,F_{q}^{1}}(\frak{b}^{(1)})$ we have two
group like elements
\begin{equation*}
\begin{array}{c}
Z_{1}=Wq^{2h_{23}^{\perp }}W^{-1}=q^{2h_{23}^{\perp }}\frac{1}{%
1-(q^{-2}-1)s^{-2}t~e_{13}} \\[2ex]
Z_{2}=Wq^{2h_{13}^{\perp }}W^{-1}=q^{2h_{13}^{\perp }}\frac{1}{%
1+(q^{2}-1)~s^{-1}t~q^{-h_{23}^{\perp }}e_{23}}.
\end{array}
\end{equation*}
This allows us to use also the Abelian twist $q^{s^{-1}\ln (Z_{1})\otimes
\ln (Z_{2})}$ . The product
\begin{equation*}
F_{q}=q^{s^{-1}\ln (Z_{1})\otimes \ln (Z_{2})}\cdot F_{q}^{1}
\end{equation*}
defines a $q-$twist with the property
\begin{equation}
F_{q}\equiv \exp ((h_{23}^{\perp }+e_{23})\otimes \ln
(1+2t^{2}e_{13}^{\prime }))~{\rm mod}\left(st \right)  \label{qproperty}
\end{equation}

Notice that in the limit $q\longrightarrow 1$ the quantum twist structures
for $\frak{b}_{\lambda }$, $\frak{b}^{(1)}$ and $\frak{b}^{(0)}$
degenerate. They lead to equivalent families of ordinary Jordanian twists (%
\ref{jord}) $\mathcal{F}_{\mathcal{J}}=\exp (H\otimes \sigma \left( \xi
\right) ,\quad \sigma \left( \xi \right) =~\ln (1+~\xi E).$

\subsection{Quantum twists with four-dimensional carriers}

Similar to the previous study (Section 2) we consider separately the
nonequivalent classes of four-dimensional Lie Frobenius subalgebras:
\begin{equation*}
\frak{r}=\left( {\
\begin{array}{rrr}
\ast  & \ast  & \ast  \\[2ex]
0 & \ast  & 0 \\[2ex]
0 & 0 & \ast
\end{array}
}\right)
\end{equation*}
and
\begin{equation*}
\frak{q}_{\,a_{1},a_{2},a_{3}}=\ast \left( {\
\begin{array}{rrr}
a_{1} & 0 & 0 \\[2ex]
0 & a_{2} & 0 \\[2ex]
0 & 0 & a_{3}
\end{array}
}\right) +\left( {\
\begin{array}{rrr}
0 & \ast  & \ast  \\[2ex]
0 & 0 & \ast  \\[2ex]
0 & 0 & 0
\end{array}
}\right)
\end{equation*}
We had one family of solutions associated to $\frak{r}$ and three
nonequivalent classes associated with a particular choice of $%
(a_{1},~a_{2},~a_{3})$ (see (\ref{norm-ext}),(\ref{deff-jord}),(\ref{per-ext}%
) and (\ref{two-cart-r})).

\subsubsection{Case $\frak{r}$}

Due to the isomorphism
\begin{equation*}
\frak{r}\cong \left( {\
\begin{array}{rrr}
\ast  & 0 & \ast  \\[2ex]
0 & \ast  & \ast  \\[2ex]
0 & 0 & \ast
\end{array}
}\right)
\end{equation*}
the case $\frak{r}$ can be treated similarly to $\frak{b^{(1)}}$ . Define
the $q-$twist as the coboundary twist
\begin{eqnarray*}
F_{q} &=&(W\otimes W)\Delta (W^{-1}), \\
W &=&\exp _{q^{-2}}(-s^{-1}t~q^{-h_{23}^{\perp }}e_{23})\exp
_{q^{2}}(s^{-1}t~e_{13})
\end{eqnarray*}
and assume that $q\equiv 1+st~{\rm mod}\left(s^{2}t \right)$. Explicitly,
\begin{equation*}
\begin{array}{r}
F_{q}=\exp _{q^{-2}}(-s^{-1}t~1\otimes q^{-h_{23}^{\perp }}e_{23})\exp
_{q^{2}}(s^{-1}t~q^{-2h_{13}^{\perp }}\otimes q^{-h_{23}^{\perp
}}e_{23})\cdot  \\[2ex]
\cdot \exp _{q^{2}}(s^{-1}t~1\otimes e_{13}^{\prime })\exp
_{q^{-2}}(-s^{-1}t~q^{-2h_{23}^{\perp }}\otimes e_{13}^{\prime })
\end{array}
\end{equation*}
and
\begin{equation*}
F_{q}\equiv \exp (-h_{13}^{\perp }\otimes \ln (1+2t^{2}~e_{23}))\exp
(h_{23}^{\perp }\otimes \ln (1+2t^{2}~e_{13}))~{\rm mod}\left(st \right).
\end{equation*}
In the limit $s\longrightarrow 0$ this expression gives the double-Jordanian
twist (\ref{two-cart-r}).

\subsubsection{Case $\frak{q}_{a_{1},a_{2},a_{3}}$}

It was shown in the table (\ref{cohomol}) that we can subdivide the case $\frak{q}%
_{\,a_{1},a_{2},a_{3}}$ into the subclasses according to their cohomological
properties:
\begin{equation*}
H^{2}(\frak{q}_{\,a_{1},a_{2},a_{3}})=\left\{
\begin{array}{l}
(a_{1},a_{2},a_{3})=(0,1,-1),(1,1,-2) \\[2ex]
0\ \mathrm{otherwise}
\end{array}
\right.
\end{equation*}
Apply the Abelian twist
\begin{equation*}
{\mathcal K}=q^{-(2\zeta +1)~h_{13}^{\perp }\otimes h_{23}^{\perp }}
\end{equation*}
to the Hopf algebra $U_{q}(\frak{sl}_{3})$. Define $U_{q}(\frak{q}%
_{a_{1},a_{2},a_{3}})$ as the minimal Hopf subalgebra in $U_{q,{\mathcal K}}(\frak{sl}%
_{3})$ containing $e_{13},e_{12},e_{23}$ and $h_{\zeta }=h_{23}^{\perp
}+\zeta ~h_{13}^{\perp }$. The coproduct of $e_{13}$ in $U_{q,{\mathcal K}}(\frak{sl}%
_{3})$ has the form
\begin{equation*}
\Delta (e_{13})=q^{-2h_{\zeta }}\otimes e_{13}+e_{13}\otimes
1+(1-q^{2})~e_{12}q^{-h_{23}}\otimes q^{-(2\zeta +1)~h_{23}^{\perp }}e_{23}.
\end{equation*}
The $q-$twist is defined by the coboundary expression
\begin{equation*}
F_{q}=(W\otimes W)\Delta (W^{-1}),\qquad W=\exp _{q^{2}}(ts^{-1}~e_{13})
\end{equation*}
and the limit $s\rightarrow 0$ taken along the curve
\begin{equation*}
q\equiv 1+st~{\rm mod}\left(s^{2}t \right).
\end{equation*}
gives two types of twists (\ref{norm-ext}) and (\ref{per-ext}).

Now consider $\frak{q}_{\,0,1,-1}$. The corresponding $r-$matrix has the
following form:
\begin{equation*}
r_{1}(\eta )=H_{23}\wedge E_{23}+2\eta ~E_{12}\wedge E_{13}.
\end{equation*}
It is equivalent to the $r-$matrix
\begin{equation*}
r_{2}(\lambda )=H_{23}\wedge E_{23}+\lambda ~(H_{23}\wedge
E_{13}+E_{12}\wedge E_{23})
\end{equation*}
via the transformation:
\begin{equation*}
r_{2}(i\sqrt{2\eta })=\mathrm{Ad}\exp (i\sqrt{2\eta }~E_{12})\otimes \exp (i%
\sqrt{2\eta }~E_{12})\circ(r_{1}(\eta )).
\end{equation*}
We can propose that $U_{q}(\frak{q}_{\,0,1,-1})$ is just a Hopf subalgebra in $%
U_{q}(\frak{sl}_{3})$ spanned by $\{e_{12},e_{23},q^{\pm h_{12}},q^{\pm
h_{23}}\}$. Though it seems that there is no easy way to obtain $F_{q}$ that
contain the factors
\begin{equation*}
\exp (\frac{1}{2}H_{23}\otimes \ln (1+tE_{23}-\frac{1}{2}\eta
^{2}t^{3}~E_{13}^{2}))\exp (\eta ^{2}t^{2}~E_{12}\otimes E_{13}).
\end{equation*}
necessary to guarantee the desired properties.

\subsection{Quantum twist with six-dimensional carrier}

As it was mentioned above up to the conjugation the only six-dimensional
subalgebra is
\begin{equation*}
\frak{p}=\left( {\
\begin{array}{rrr}
\ast  & \ast  & \ast  \\[2ex]
0 & \ast  & \ast  \\[2ex]
0 & \ast  & \ast
\end{array}
}\right)
\end{equation*}
Let $U_{q}(\frak{p})$ be a Hopf subalgebra in the algebra $U_{q,{\mathcal K}}(\frak{sl%
}_{3})$ obtained as a deformation of $U_{q}(\frak{sl}_{3})$ by the Abelian
twist $K$,
\begin{equation*}
{\mathcal K}=q^{h_{12}^{\perp }\otimes h_{23}^{\perp }}
\end{equation*}
In the subalgebra $U_{q}(\frak{p})$ we have the following coproducts:
\begin{equation*}
\begin{array}{l}
\Delta _{\mathcal K}(e_{12})=q^{-2h_{13}^{\perp }}\otimes e_{12}+e_{12}\otimes 1, \\%
[2ex]
\Delta _{\mathcal K}(q^{-h_{23}^{\perp }}e_{23})=q^{-2h_{12}^{\perp }}\otimes
q^{-h_{23}^{\perp }}e_{23}+q^{-h_{23}^{\perp }}e_{23}\otimes 1, \\
\Delta _{K}(e_{32})=e_{32}\otimes q^{-2h_{13}^{\perp }}+1\otimes e_{32}.
\end{array}
\end{equation*}
It follows, \cite{KM}, that the element
\begin{equation*}
F_{{\mathcal K}M}=\exp _{q^{-2}}((q-q^{-1})t~e_{32}\otimes e_{12}).
\end{equation*}
 is a twist for $U_{q}(\frak{p})$.

Let us consider the twist $F_{q}$ equivalent to $F_{{\mathcal K}M}$ ,
\begin{equation*}
F_{q}=(W\otimes W)F_{{\mathcal K}M}\Delta (W^{-1})
\end{equation*}
here
\begin{equation*}
W=\exp _{q^{2}}(t^{2}s^{-1}~e_{12})\exp _{q^{2}}(ts^{-1}~q^{-h_{23}^{\perp
}}e_{23})\exp _{q^{2}}(t^{2}s^{-1}~e_{12}).
\end{equation*}
Explicitly we have
\begin{equation*}
\begin{array}{r}
F_{q}=(\exp _{q^{2}}(t^{2}s^{-1}~e_{12})\exp
_{q^{2}}(ts^{-1}~q^{-h_{23}^{\perp }}e_{23})\otimes W)\cdot  \\[2ex]
\cdot \exp _{q^{-2}}((q-q^{-1})t~e_{32}\otimes e_{12})\cdot \exp
_{q^{-2}}(-t^{2}s^{-1}~q^{-2h_{13}^{\perp }}\otimes e_{12})\cdot  \\[2ex]
\cdot \exp _{q^{-2}}(-ts^{-1}~q^{-h_{23}}e_{23}\otimes 1)\cdot \exp
_{q^{-2}}(-ts^{-1}~q^{-2h_{12}^{\perp }}\otimes q^{-h_{23}^{\perp
}}e_{23})\cdot  \\
\cdot \exp _{q^{-2}}(-t^{2}s^{-1}~e_{12}\otimes 1)\exp
_{q^{-2}}(-t^{2}s^{-1}~q^{-2h_{13}^{\perp }}\otimes e_{12}).
\end{array}
\end{equation*}
To transform $F_{q}$ further we use the commutation property
\begin{equation*}
\lbrack q^{-h_{23}^{\perp }}e_{23},e_{32}]=\frac{q^{-2h_{13}^{\perp
}}-q^{-2h_{12}^{\perp }}}{q-q^{-1}}
\end{equation*}
and the relations
\begin{equation*}
\begin{array}{l}
\exp _{q^{-2}}((q-q^{-1})t~e_{32}\otimes e_{12})\exp
_{q^{-2}}(-ts^{-1}~q^{-2h_{12}^{\perp }}\otimes q^{-h_{23}^{\perp }}e_{23})=
\\[2ex]
\exp _{q^{-2}}(-ts^{-1}~q^{-2h_{12}^{\perp }}\otimes q^{-h_{23}^{\perp
}}e_{23})\exp _{q^{-2}}((1-q^{2})s^{-1}t^{2}~e_{32}q^{-2h_{12}^{\perp
}}\otimes q^{-h_{23}^{\perp }}e_{13}^{\prime })\cdot  \\[2ex]
\cdot \exp _{q^{-2}}((q-q^{-1})t~e_{32}\otimes e_{12}).
\end{array}
\end{equation*}
As a result the corresponding factors in $F_{q}$ can be transposed,
\begin{equation*}
\begin{array}{r}
\exp _{q^{-2}}((q-q^{-1})t~e_{32}\otimes e_{12})\exp
_{q^{-2}}(-t^{2}s^{-1}~q^{-2h_{13}^{\perp }}\otimes e_{12}) \times
\rule{3mm}{0mm} \\
\times \exp
_{q^{-2}}(-ts^{-1}~q^{-h_{23}^{\perp }}e_{23}\otimes 1)= \\[2ex]
=\exp _{q^{-2}}(-ts^{-1}~q^{-h_{23}}e_{23}\otimes 1)\exp
_{q^{-2}}(-t^{2}s^{-1}~q^{-2h_{12}^{\perp }}\otimes e_{12})\times
\rule{3mm}{0mm}\\
\times
\exp
_{q^{-2}}((q-q^{-1})t~e_{32}\otimes e_{12}),
\end{array}
\end{equation*}
and the twisting element takes the form
\begin{equation*}
\begin{array}{l}
F_{q}=(1\otimes W)\exp _{q^{-2}}(-t^{2}s^{-1}~q^{-2h_{12}^{\perp }}\otimes
e_{12})\times
\\
\rule{6mm}{0mm}\times\exp _{q^{-2}}(-ts^{-1}~q^{-2h_{12}^{\perp }}\otimes
q^{-h_{23}^{\perp }}e_{23})\times
\\
\rule{6mm}{0mm}\times
\exp _{q^{-2}}((1-q^{2})s^{-1}t^{2}~e_{32}q^{-2h_{12}^{\perp }}\otimes
q^{-h_{23}^{\perp }}e_{13}^{\prime })\times
\\
\rule{6mm}{0mm}\times\exp
_{q^{-2}}((q-q^{-1})t~e_{32}\otimes e_{12})
 \exp _{q^{-2}}(-t^{2}s^{-1}~q^{-2h_{13}^{\perp }}\otimes e_{12}).
\end{array}
\end{equation*}
Now taking into account that the following commutator iz zero,
\begin{equation*}
\begin{array}{l}
\lbrack \exp _{q^{2}}(ts^{-1}~1\otimes q^{-h_{23}}e_{23})\exp
_{q^{2}}(t^{2}s^{-1}~1\otimes e_{12}), \\[2ex]
\exp _{q^{-2}}(-t^{2}s^{-1}~q^{-2h_{12}^{\perp }}\otimes e_{12})\exp
_{q^{-2}}(-ts^{-1}~q^{-2h_{12}^{\perp }}\otimes q^{-h_{23}^{\perp
}}e_{23})]=0,
\end{array}
\end{equation*}
the final expression for the q-twisting element can be obtained:
\begin{equation*}
\begin{array}{l}
F_{q}=\exp _{q^{2}}(t^{2}s^{-1}~1\otimes e_{12})\exp
_{q^{-2}}(-t^{2}s^{-1}~q^{-2h_{12}^{\perp }}\otimes e_{12})\times
\\
\rule{3mm}{0mm}\times
 \exp _{q^{2}}(ts^{-1}~1\otimes q^{-h_{23}^{\perp }}e_{23})\cdot \exp
_{q^{-2}}(-ts^{-1}~q^{-2h_{12}^{\perp }}\otimes q^{-h_{23}^{\perp
}}e_{23})\times
\\
\rule{3mm}{0mm}\times
 \exp _{q^{-2}}((1-q^{2})s^{-1}t^{2}~e_{32}q^{-2h_{23}^{\perp }}\otimes
q^{-h_{23}^{\perp }}e_{13}^{\prime })\times
\\
\rule{3mm}{0mm}\times \exp
_{q^{-2}}((q-q^{-1})t~e_{32}\otimes e_{12})\times
\\
\rule{3mm}{0mm}\times
 \exp _{q^{2}}(t^{2}s^{-1}~1\otimes e_{12})\exp
_{q^{-2}}(-t^{2}s^{-1}~q^{-2h_{13}^{\perp }}\otimes e_{12}).
\end{array}
\end{equation*}
Assuming $q\equiv 1+st~{\rm mod}\left(s^{2}t \right)$ and applying the Heine's formula
we can calculate the limit $s\rightarrow 0$ which gives the twist:
\begin{equation*}
\begin{array}{l}
\exp (H_{12}^{\perp }\otimes \ln (1+2t^{3}~E_{12}))\exp (H_{12}^{\perp
}\otimes \ln (1+2t^{2}~E_{23}))\times  \\
\times \exp (-2t^{3}E_{32}\otimes E_{13})\exp (H_{13}^{\perp }\otimes \ln
(1+2t^{3}E_{12})).
\end{array}
\end{equation*}
This expression is the special case of the parabolic twist obtained in \cite{LS}
and presented here in Section 2 (see (\ref{par-param})).

\section{Conclusions}
We have shown that the factorization property presents the possibility to
obtain all the solutions to the twist equations for algebra $U(sl(3))$.
The full list of the antisymmetric classical $r$-matrices, constant solutions of
CYBE, was quantized and the corresponding twists were constructed explicitly
in the form of product of twisting factors. Each of these factors refers to an
independent solution of the twist equation.

We have also demonstrated that when the Drinfeld-Jimbo $r$-matrix and the
antisymmetric $r$-matrix corresponding to the twist are compatible (that is
their sum gives rise to a solution for the modified classical Yang-Baxter
equation) the quantum counterpart of this twist can be obtained.

It is known that triangular twists permit to deform
integrable models related to Yangians \cite{KSt}.
We suppose that constructed coboundary $q$-analogues of triangular twists
give rise to a possibility to study mentioned above deformed models
starting directly from known anisotropic models. The latter being connected
with the corresponding quantum affine algebras are similarly transformed
under coboundary twists. Hence a new basis of eigenvectors will appear.

\section{Acknowledgements}

This work was supported in part by the RFBR grant N 06-01-00451,
the CRDF grant No RUMI-2622-ST-04 and the program RNP 2.1.1 grant No 1112.

\appendix
\section*{Appendix}
\setcounter{section}{1}

Contrary to the situation described in the beginning of the subsection 2.6.1
the peripheric carriers $\mathbf{L}_{1,0}$ ($\mathbf{L}_{0,1}$) present more
interesting possibilities -- the twisting elements can be enlarged by the
additional factors. These constructions can be proved to be equivalent to
the ordinary double-Cartan case described in 1.3.2 but deserve separate
presentation. Here the carrier $\mathbf{L}_{1,0}$ is more convenient for our
purposes and the peripheric twist looks like
\begin{equation}
\mathcal{F}_{\mathcal{P}}=\exp (\xi E_{12}\otimes E_{23})\exp (H_{\mathcal{P}%
}\otimes \sigma (\xi )),  \label{per-ext22}
\end{equation}
with
\begin{equation}
H_{\mathcal{P}}=\frac{2}{3}E_{11}-\frac{1}{3}E_{22}-\frac{1}{3}E_{33},\quad
\quad \sigma (\xi )=\ln (1+\xi E_{13})
\end{equation}
The costructure can be obtained from (\ref{delf}):
\begin{equation}
\begin{array}{lcl}
\Delta _{\mathcal{P}}(H_{\mathcal{P}}) & = & H_{\mathcal{P}}\otimes
e^{-\sigma }+1\otimes H_{\mathcal{P}}-E_{12}\otimes E_{23}e^{-\sigma }, \\%
[0.2cm]
\Delta _{\mathcal{P}}(E_{12}) & = & E_{12}\otimes 1+1\otimes E_{12}, \\%
[0.2cm]
\Delta _{\mathcal{P}}(E_{23}) & = & E_{23}\otimes +e^{\sigma }\otimes E_{23},
\\[0.2cm]
\Delta _{\mathcal{P}}(E_{13}) & = & E_{13}\otimes e^{\sigma }+1\otimes
E_{13}.
\end{array}
\label{delp2}
\end{equation}
The element $E_{12}$ remains primitive. Now in $\frak{g}$ there is a Cartan
element
\begin{equation}
H^{\perp }=\frac{1}{3}E_{11}-\frac{2}{3}E_{22}+\frac{1}{3}E_{33},
\end{equation}
whose dual is orthogonal to $\lambda =e_{1}-e_{3}$. Consequently this
element also remains primitive after the peripheric twist $\mathcal{F}_{%
\mathcal{P}}$. The corresponding Borel subalgebra $\mathbf{B}(2)$ with the
generators $H^{\perp },E_{12}$ can be twisted additionally by the Jordanian
twist (\ref{jord}). As a result the triple of twisting factors form a twist
\begin{equation}
\mathcal{F}_{\mathcal{JPJ}}=\exp (H^{\perp }\otimes \sigma _{12}(\xi ))\exp
(\xi E_{12}\otimes E_{23})\exp (H_{\mathcal{P}}\otimes \sigma (\xi )),
\label{add-jord}
\end{equation}
with
\begin{equation}
\sigma _{12}(\xi )=\ln (1+\xi E_{12}).
\end{equation}
Still the carrier algebra for this twist cannot have the dimension greater
than four. In the classical $r$-matrix the additional term originating from
the factor $\exp (\xi H^{\perp }\otimes \sigma _{12})$ induces a change of
the $E_{23}$ basic element for $B=E_{23}-\xi H^{\perp }$. With this change
the four-dimensional space of the carrier generated by $\left\{ H_{\mathcal{P%
}},E_{12},E_{13},B=H^{\perp }-\frac{1}{\xi }E_{23}\right\} $ becomes closed
under the compositions of $\frak{g}$. Thus we obtain the deformation $%
\mathbf{L}_{1,0}^{\mathrm{def}}$ with the relations
\begin{equation}
\begin{array}{ll}
\left[ H_{\mathcal{P}},E_{12}\right] =E_{12}, & \left[ H_{\mathcal{P}},B%
\right] =0, \\
\left[ H_{\mathcal{P}},E_{13}\right] =E_{13}, & \left[ E_{12},E_{13}\right]
=0, \\
\left[ E_{12},B\right] =E_{13}+\xi E_{12}. &
\end{array}
\label{def-carrier}
\end{equation}
Despite the fact that the deforming function $\mu $ with $\mu
(E_{12},B)=E_{12}$ is a coboundary ($\mu \in B^{2}(\mathbf{L}_{1,0},\mathbf{L%
}_{1,0})$) this deformation is nontrivial. This can be checked by inspecting
the ranks: $\mathrm{rank}(\mathbf{L}_{1,0})=1$ and $\mathrm{rank}(\mathbf{L}%
_{1,0}^{\mathrm{def}})=2$. (Notice that the similarity transformation that
cancels $\mu $ brings the cohomologically nontrivial term in the second
order of the deformation parameter.) Thus incorporating the generator $%
H^{\perp }$ in the structure of the twist we have passed to the new carrier $%
\mathbf{L}_{1,0}^{\mathrm{def}}$. The latter must be identified with a
Frobenius subalgebra in $\frak{g}$. To find such consider the new basis $%
\left\{ H=H_{\mathcal{P}}-B,B,A=\xi E_{12}+E_{13},E=E_{13}\right\} $. Now
the commutation relations are
\begin{equation}
\begin{array}{lll}
\left[ H,A\right] =0, & \left[ B,A\right] =A, & \left[ H,B\right] =0, \\
\left[ H,E\right] =E, & \left[ B,E\right] =0, & \left[ A,E\right] =0,
\end{array}
\label{def-carrier2}
\end{equation}
Thus $\mathbf{L}_{1,0}^{\mathrm{def}}=\mathbf{B}(2)\oplus \mathbf{B}(2)$
(the structure that we had in (\ref{two-bor})). Obviously having this form
for $\mathbf{L}_{1,0}^{\mathrm{def}}$ we can apply to it the double-Jordanian
twist (\ref{two-cart}):
\begin{equation}
\mathcal{F}_{\mathcal{JJ}}=\exp (H\otimes \sigma _{E})\exp (B\otimes \sigma
_{A})\quad \quad \sigma _{A}=\ln (1+A),\sigma _{E}=\ln (1+E).
\label{two-cart2}
\end{equation}
Returning to the initial basis in $\frak{g}$ it can be written as
\begin{equation}
\begin{array}{lcl}
\mathcal{F}_{\mathcal{JJ}} & = & \exp (H^{\perp }\otimes (\sigma _{(\xi
E_{12}+E_{13})}-\sigma _{13})-\frac{1}{\xi }E_{23}\otimes (\sigma _{(\xi
E_{12}+E_{13})}-\sigma _{13})) \\
&  & \exp (H_{\mathcal{P}}\otimes \sigma _{13}) \\
& = & \exp (H^{\perp }\otimes (\sigma _{(\xi E_{12}+E_{13})}-\sigma
_{13}))\exp (-E_{23}\otimes E_{12}e^{-\sigma _{13}}) \\
&  & \exp (H_{\mathcal{P}}\otimes \sigma _{13}) \\
& = & \exp (H^{\perp }\otimes \sigma _{(\xi E_{12}+E_{13})})\exp
(-E_{23}\otimes E_{12}) \\
&  & \exp ((H_{\mathcal{P}}-H^{\perp })\otimes \sigma _{13})
\end{array}
\label{two-cart3}
\end{equation}
The form of the twisting element looks similar to that of (\ref{add-jord})
but is different. Notice that here the first two twisting factors present
the peripheric twist that produces the primitive coproduct for the element $%
E_{23}$. So the last quasi-Jordanian factor is based on the quasiprimitive
combination of elements. The $r$-matrix is the same as in the case (\ref
{add-jord}).

\end{document}